\documentclass[12pt]{amsart}
\usepackage{fancyhdr} 
\usepackage{tikz}
\usetikzlibrary{cd}
\usepackage{amsmath}
\usepackage[active]{srcltx}
\usepackage[pagebackref, colorlinks, linkcolor=red, citecolor=blue, urlcolor=blue, hypertexnames=true]{hyperref}
\usepackage[backrefs]{amsrefs}
\usepackage{enumerate}
\usepackage[mathscr]{eucal}
\usepackage{appendix}
\allowdisplaybreaks
\usepackage[all,cmtip]{xy}

\usepackage{color}

\allowdisplaybreaks
\usepackage[all]{xy}

\newcommand{\R}{{\mathbb R}}

\newcommand{\Z}{{\mathbb Z}}
\newcommand{\N}{{\mathbb N}}

\def\m{\mu}
\def\n{\nu}
\def \d{\delta}
\def\v{\varepsilon}

\def\sub{\subset}

\newtheorem{thm}{Theorem}[section]
\newtheorem{cor}[thm]{Corollary}
\newtheorem{lem}[thm]{Lemma}

\newtheorem{definition}[thm]{Definition}

\newtheorem{remark}[thm]{Remark}
\newtheorem{pro}[thm]{Proposition}

\newtheorem{claim}[thm]{Claim}

\theoremstyle{definition}


\begin{document}
	\title{Mean Hausdorff Dimension of Some Infinite Dimensional Fractals for Amenable Group  Actions}

\email{lixq233@mail2.sysu.edu.cn }
\email{luoxf29@mail2.sysu.edu.cn }

\maketitle
\centerline{\scshape Xianqiang Li, Xiaofang Luo}

	\maketitle
	\medskip
	\noindent{\bf Abstract} For the countable discrete amenable group actions, we calculate  the mean Hausdorff dimensions of three types of infinite dimensional fractal systems, the self-similar systems, homogeneous systems in the infinite-dimensional torus, and the infinite dimensional-carpets.

	\section{Introduction}
	\hspace{4mm}
	Fractal geometry is originally related to describe the intricate natural landscapes, such as mountains, rivers, coastlines, and lakes, which  studies complex structures and shapes. Hausdorff dimension and Minkowski dimension are  two vital  quantities  for us to 
	 quantify the ``size" of sets, and allow us to  understand the geometric structures of fractals. Among the various types of fractal objects, there are three types of finite-dimensional systems that are  known as homogeneous sets, self-similar sets, and Bedford–McMullen carpets. Furstenberg \cite{FHH} and Falconer \cite{FKJ} showed that the Hausdorff dimension and Minkowski dimension of the homogeneous sets and self-similar sets are equal. However, Bedford \cites{BT} and McMullen \cites{MC} further showed that the  Hausdorff dimension and Minkowski dimension of Bedford–McMullen carpet do not coincide. On the one hand, there is a significant relationship between dimension and entropy. For instance, for the one-sided subshift of finite type $X$ of symbol space, the Hausdorff dimension of $X$ coincides with the topological entropy of $X$ with respect to the standard metric \cites{RB,FHH}, which is also valid for the actions of $\mathbb{N}^d$ and $\mathbb{Z}^d$ in Simpson's work \cite{SS}. He also mention  whether one can obtain some more precise relations between dimension and entropy in the framework of amenable groups or more general groups actions. On the other hand, it is natural to ask whether we can calculate the aforementioned three types of fractals in the case of infinite-dimensional.

	 \hspace{4mm} 
	 To capture the dynamics of infinite entropy systems,
	Lindenstrauss and Weiss \cite{TM1} introduced the metric mean dimension $\operatorname{mdim}_{\mathrm{M}}(\cdot)$, which is analogous to the definition of Minkowski dimension. An analogue of the Hausdorff dimension, known as the mean Hausdorff dimension $\operatorname{mdim}_{\mathrm{H}}(\cdot)$, was  introduced  by Lindenstrauss and Tsukamoto \cite{LT1} to link the geometric measure theory and rate-distortion theory while obtaining the double variational principle of mean dimension. There are few results involving the fractal dimensions of infinite-dimensional sets. In \cite{TM}, Tsukamoto   gave the nice framework of the above three types of fractal sets in sense of  infinite-dimensional  case, which can be stated as follows:

		\hspace{4mm} 	 
	(1) Suppose $\mathbb{R}$ is a circle with the metric $\rho $, and the metric $d$ on $(\mathbb{R} / \mathbb{Z})^{\mathbb{N}}$ is defined by $$d\left(\left(x_{n}\right)_{n \in \mathbb{N}},\left(y_{n}\right)_{n \in \mathbb{N}}\right)=\sum_{n=1}^{\infty} 2^{-n} \rho\left(x_{n}, y_{n}\right).$$
	Given a positive integer $b$ with $b\geq 2.$ The ``$\times b$ map'' $S_b$ on $(\mathbb{R} / \mathbb{Z})^{\mathbb{N}}$ is defined by $$S_{b}\left(\left(x_{n}\right)_{n \in \mathbb{N}}\right)=\left(b x_{n}\right)_{n \in \mathbb{N}}.$$
	The shift map $\sigma$ on  $(\mathbb{R} / \mathbb{Z})^{\mathbb{N}}$ is defined by $$\sigma\left(\left(x_{n}\right)_{n \in \mathbb{N}}\right)=\left(x_{n+1}\right)_{n \in \mathbb{N}}.$$
	We call the closed subset $Y\subset (\mathbb{R} / \mathbb{Z})^{\mathbb{N}}$ a \textbf{homogeneous system} if  $S_bY\subset Y$ and $\sigma Y\subset Y$.\\
	
		\hspace{4mm} 	
	(2) Suppose $c $ is a positive number with  $c<1$, $(\Omega, T)$ is a dynamical system ($\Omega$ is compact and $T$ is continous map) and $\ell^{\infty}$ is the  space of bounded sequences with norm $\|x\|_{\infty}:=\sup _{n \geq 1}\left|x_{n}\right|$  for  $x=\left(x_{n}\right)_{n \in \mathbb{N}}\in \ell^{\infty}.$ The shift map  $\sigma$ on $\ell^{\infty}$ is continuous with respect to weak$^{*}$ topology. For any $\omega \in \Omega$, there is a point $a(\omega)=\left(a(\omega)_{n}\right)_{n \in \mathbb{N}} \in \ell^{\infty}$  such that the map
	$$\Omega \ni \omega \mapsto a(\omega) \in \ell^{\infty}$$
	is continuous (with respect to the weak${ }^{*} $ topology of  $\ell^{\infty} $ ) and equivariant (i.e.  $\sigma(a(\omega))=   a(T \omega)$). And for any  $\omega \in \Omega$, there is a contracting similarity transformation  $S_{\omega}: \ell^{\infty} \rightarrow \ell^{\infty}$ which defined by
	$$S_{\omega}(x)=c x+a(\omega).$$ We call the non-empty compact subset  $X$  of  $\ell^{\infty}$ the \textbf{self-similar system} if $X$ is the unique set satisfying $$X=\bigcup_{\omega \in \Omega} S_{\omega}(X).$$\\
	
		\hspace{4mm} 	
	(3) Suppose $a,b$ are two natural numbers with $a \geq b \geq 2$. Setting
	$$A=\{0,1,2, \ldots, a-1\}, \quad B=\{0,1,2, \ldots, b-1\}.$$
	Let  $(A \times B)^{\mathbb{N}}$  be the one-sided full-shift on the alphabet  $A \times B$  with the shift map  $\sigma:(A \times B)^{\mathbb{N}} \rightarrow(A \times B)^{\mathbb{N}}$, $\pi:(A \times B)^{\mathbb{N}} \rightarrow B^{\mathbb{N}}$  be the natural projection and $[0,1]^{\mathbb{N}}=[0,1] \times[0,1] \times[0,1] \times \ldots$ be the infinite dimensional cube. The metric $d$ of the product  $[0,1]^{\mathbb{N}} \times[0,1]^{\mathbb{N}} $ is defined by
	$$d\left((x, y),\left(x^{\prime}, y^{\prime}\right)\right)=\sum_{n=1}^{\infty} 2^{-n} \max \left(\left|x_{n}-x_{n}^{\prime}\right|,\left|y_{n}-y_{n}^{\prime}\right|\right),$$
	where  $x=\left(x_{n}\right)_{n \in \mathbb{N}}, y=\left(y_{n}\right)_{n \in \mathbb{N}} $ and  $x^{\prime}=\left(x_{n}^{\prime}\right)_{n \in \mathbb{N}}, y^{\prime}=\left(y_{n}^{\prime}\right)_{n \in \mathbb{N}}$  are points in  $[0,1]^{\mathbb{N}} .$ The shift map $\sigma$ on  $[0,1]^{\mathbb{N}}\times[0,1]^{\mathbb{N}} $ is defined by
	$$\sigma\left(\left(x_{n}\right)_{n \in \mathbb{N}},\left(y_{n}\right)_{n \in \mathbb{N}}\right)=\left(\left(x_{n+1}\right)_{n \in \mathbb{N}},\left(y_{n+1}\right)_{n \in \mathbb{N}}\right)$$
	If $\emptyset\neq \Omega \subset(A \times B)^{\mathbb{N}}$ is a subshift, then we call $X_{\Omega} \subset[0,1]^{\mathbb{N}} \times[0,1]^{\mathbb{N}}$:  
	$$X_{\Omega}=\left\{\left.\left(\sum_{m=1}^{\infty} \frac{x_{m}}{a^{m}}, \sum_{m=1}^{\infty} \frac{y_{m}}{b^{m}}\right) \in[0,1]^{\mathbb{N}} \times[0,1]^{\mathbb{N}} \right\rvert\,\left(x_{m}, y_{m}\right) \in \Omega \text { for all } m \geq 1\right\}$$
	a \textbf{carpet system}.

		\hspace{4mm}  	
  Tsukamoto \cite{TM} calculated the mean Hausdorff dimension  and metric mean dimension of these infinite-dimensional fractal set, and obtained the  similar results as the finite-dimensional case.   More work  toward this  aspect, see  \cites{LF, TM, LSL1, LSL2} for the extension of  Tsukamoto's work by using the preimage metric mean dimension and   mean $\psi$-intermediate dimension. In this paper,   we are interested in the calculation of  the mean  Hausdorff dimension of  ``infinite-dimensional fractals'' in the context of  amenable groups actions.

	 \hspace{4mm}
	 	Let $G$ be a countable discrete amenable group. By a pair $(X,G,\Gamma)$ we mean a $G$-system, where $X$ is a compact metric space and $\Gamma:G\times\tilde{X}\to X$, given by $(g,x)\to$ $gx$, is a continuous mapping satisfying:\\
	 $(1)\Gamma(1_{G},x)=x$ for every $x\in X;$\\
	 $(2)\Gamma(g_{1},\Gamma(g_{2},x))=\Gamma(g_{1}g_{2},x)$ for every $g_1,g_{2}\in G$ and $x\in X.$\\
	 When there is no ambiguity on the map $\Gamma$, we write $(X,G)$ instead of $(X,G,\Gamma)$.

	\hspace{4mm}
	The main result of  this paper can be stated as follows:

	\begin{thm}\label{t1.1}
		For the subsystem $\left(X_{\Omega}, \sigma\right)$ (carpet system) of $\left([0,1]^{G} \times[0,1]^{G}, \sigma\right),$ we have 
		\begin{align*}
			&\operatorname{mdim}_{\mathrm{H}}\left(X_{\Omega}, \{B_S(m)\},  d\right) =\frac{h_{\text {top }}^{w}(\pi, G)}{\log b},\\
			&\operatorname{mdim}_{\mathrm{M}}\left(X_{\Omega}, G, d\right)= \frac{h_{\text {top }}(\Omega, G)}{\log a}+\left(\frac{1}{\log b}-\frac{1}{\log a}\right) h_{\text {top }}\left(\Omega^{\prime}, G\right).
		\end{align*}
		Here $h_{\text {top }}(\cdot)$ and $h_{\text {top }}^{w}(\cdot)$ denote the topological entropy and the weighted topological entropy.
	\end{thm}

	\begin{thm}\label{t1.2}
		Let $G$ be a countable discrete amenable group and $( \Omega,G, T) $ be a $G$-system. For the self-similar dynamical system $(X,G,\sigma)$ defined by a family of contractions $\{S_{\omega}\}_{\omega \in \Omega}$ in $\ell^{\infty}$, where $$\ell^\infty=\left\{(x_g)_{g\in G}\in \mathbb{R}^{G}  \left|\sup\limits_{g\in G}|x_g|<\infty\right\}\right.$$ 
		is the space of bounded sequences with the norm $\|x\|_\infty:= \sup _{g\in G}|x_g|$ for $x= ( x_g) _{g \in  G },$ and $\ell^{\infty}$ is endowed with the weak* topology. We defined a metric d on X by
		$$d(\left(x_{g}\right)_{g \in G},\left(y_{g}\right)_{g \in G})=\sum_{g\in G} \alpha_{g} |x_{g}- y_{g}| ,$$
		where
		$\alpha_{g} \in (0,+\infty)$ satisfies $$
		\alpha_{1_{G}}=1,\sum_{g \in G}\alpha_{g} < +\infty.
		$$
		
		Let $\{F_{n}\}^{\infty}_{n=1}$ be a F$\phi$lner sequence in $G$, then
		$$
		\mathrm{mdim}_{\mathrm{H}}(X,\{F_{n}\},d)=\mathrm{mdim}_{\mathrm{M}}(X,G,d)\leq\frac{h_{\mathrm{top}}(\Omega,G)}{\mathrm{log}(1/c)}.
		$$
		Here $c$ is the ratio of contracting similarity transformation $S_{\omega}$. 
	\end{thm}
	
	\begin{thm}\label{t1.3}
		Let $\{F_{n}\}^{\infty}_{n=1}$ be a F$\phi$lner sequence in $G$. If there is a closed subset $X\subset(\mathbb{R} / \mathbb{Z})^{G}$ satisfies that $T^{m}_b(X) \subset X$ for any natural number $m$ and $ \sigma^{h}(X) \subset X$ for any $h\in G$, then we have
		$$\operatorname{mdim}_{\mathrm{H}}(X, \{F_{n}\},  d)=\operatorname{mdim}_{\mathrm{M}}(X,G,  d)=\frac{h_{\mathrm{top }}\left(X,G\times \N, \sigma, T_b\right)}{\log b} .$$
		Here  $h_{\mathrm{top }}\left(X, G \times \N, \sigma, T_b\right)$ is the topological entropy of the $G \times \mathbb{N}$-actions $\left(X, G \times \N, \sigma, T_b\right)$.
		
	\end{thm}
	
	\hspace{4mm}
	The paper is organized as follows: In Section 2, we introduce the mean dimension quantities for amenable group actions and review some basic properties. In Section 3, we compute the mean Hausdorff dimension and metric mean dimension of infinite-dimensional self-similar systems and homogeneous systems on the infinite-dimensional torus for countable discrete amenable group actions. In Section 4, we calculate the mean Hausdorff dimension and metric mean dimension of infinite-dimensional carpets for finitely generated amenable group actions. Finally, in the Appendix, we present two examples illustrating cases where the mean Hausdorff dimension and metric mean dimension can be equal or unequal.

	\section{Preliminaries}
	\hspace{4mm} 
	In the paper, let  $\N$ denote the set of natural numbers and $\#A$ denote the cardinality of the set $A$.
	
	\hspace{4mm}
	In this section, we will review the notions of amenable groups actions, mean dimension, metric mean dimension and  mean Hausdorff dimension, and we will prove an inequality for these dimensions for amenable group actions. Then we review the notion of weighted topological entropy. We end this section with introducing \texorpdfstring{$G \times \N$}{}-actions and topological Entropy for \texorpdfstring{$G \times \N$}{}-actions.
	
	\subsection{Amenable group actions}
	\
	\newline \indent \quad\ Let $G$ be a countable discrete infinite group and let $Fin(G)$ denote the set of all finite non-empty subsets of $G$. One says that a sequence $\{F_{n}\}^{\infty}_{n =1 }$ of non-empty finite subsets of $G$ is a \textbf{F$\phi$lner sequence} for $G$ if one has
	$$ \lim_{n \rightarrow \infty}\frac{|F_{n} \setminus gF_{n}|}{|F_{n}|}=0, \text{ for all } g \in G.$$
	
	The group $G$ is said to be \textbf{amenable} if it admits a F$\phi$lner sequence. It is well known that $\Z$ is amenable group.

	\hspace{4mm} 
	Now we review some results on countable amenable groups. 
	\begin{pro}
		Suppose that $G_{1}$ and $G_{2}$ are countable amenable groups then the group $G=G_{1} \times G_{2}$ is also amenable.
	\end{pro}
	
	\begin{lem} \cites{LT1,GM,OW}\label{l100}(Ornstein-Weiss lemma)
		Let $G$ be a countable amenable group and $\mathcal{F}=\{F_{n}\}^{\infty}_{n = 1}$ be a F$\phi$lner sequence for $G$.  Suppose that $h: \emph{Fin}(G) \rightarrow \mathbb{R}$ is a real-valued map satisfying the following conditions:\\
		(1) $h$ is subadditive, i.e., one has
		$$ h(E \cup F) \leq h(E)+h(F) ~for~ all~ E, F \in \emph{Fin}(G);$$
		(2) $h$ is G-invariant, i.e., one has
		$$h(Fg)=h(F)~for~all~g \in G~and~F \in \emph{Fin}(G).
		$$	
		Then the limit
		$$
		\lambda=\lim_{n \rightarrow \infty}\frac{h(F_{n})}{|F_{n}|}
		$$
		exists and one has $0 \leq \lambda < \infty$. Moreover, the limit $\lambda$ does not depend on the choice of the F$\phi$lner sequence $\mathcal{F}$ for $G$.
	\end{lem}

	\subsection{Mean dimension}
	\
	\newline \indent \quad\ Let $X$ be a compact metric space with a metric $d$. Given $\varepsilon>0$, the map $f:X \rightarrow Y$ is called an \textbf{$\varepsilon$-embedding} if diam$f^{-1}(y)< \varepsilon$ for every $y \in Y$. We define the \textbf{$\varepsilon$-width dimension} $\mathrm{Widim}_{\varepsilon}(X,d)$ as the minimum $n\geq 0$ such that there exists an $\varepsilon$-embedding $f:X \rightarrow P$ from $X$ to some $n$-dimensional simplicial complex $P$. Recall that the \textbf{topological dimension} of $X$ is given by
	$$\mathrm{dim}(X)=\lim_{\varepsilon \rightarrow 0}\mathrm{Widim}_{\varepsilon}(X,d).$$

	\hspace{4mm}
	Let $(X,G)$ be a $G$-system with a metric $d$. For $F \in \text{Fin}(G)$ and $x,y \in X$, setting
	$$ d_{F}(x,y)=\max_{g \in F}d(gx,gy).$$
	The \textbf{mean dimension} of $(X,G)$ is defined by
	$$\mathrm{mdim}(X,G)=\lim_{\varepsilon \rightarrow 0}\left( \lim_{n \rightarrow \infty}\frac{\mathrm{Widim}_{\varepsilon}(X,d_{F_{n}})}{|F_{n}|}\right).
	$$
	where $\{F_{n}\}^{\infty}_{n=1}$ is a F$\phi$lner sequence of $G$. It's easy to see that $h(F):=\text{Widim}_{\varepsilon}(X,F)$ satisfies Lemma \ref{l100}, hence the limit always exists and the value $\mathrm{mdim}(X, G)$ is independent of the choice of the F$\phi$lner sequence of $G$.
	
	\subsection{Metric mean dimension}
	\
	\newline \indent \quad\ For any positive number $\varepsilon$, the \textbf{$\varepsilon$-covering number}  $\#(X,d,\varepsilon)$ is defined by
	$$\#(X,d,\varepsilon)=\left.\min\left\{n \geq 1 \begin{array}{c|c}&\exists\text{ open covers } \{U_1,\ldots,U_n\} \text{ of } X \text{ satisfying }\\&\text{ diam}\left(U_k,d\right)<\varepsilon\quad\text{for all }1\leq k\leq n\end{array}\right.\right\}.$$
  We use the convention $\#  \left(\emptyset, d_{F_{n}}, \varepsilon\right)=0$. The \textbf{$\varepsilon$-scale Minkowski dimension} of $X$ is defined by
	$$\text{dim} _ {\text{M}} (X,d,\varepsilon)=\frac{\log\#(X,d,\varepsilon)}{\log(1 /\varepsilon)}.$$

	The \textbf{upper and lower Minkowski dimensions} of $(X,d)$ are defined by
	$$\overline{\dim}_{\mathrm{M}}(X,d)=\limsup_{\varepsilon\to0}\dim_{\mathrm{M}}(X,d,\varepsilon),~ \underline{\dim}_{\mathrm{M}}(X,d)=\liminf_{\varepsilon\to0}\dim_{\mathrm{M}}(X,d,\varepsilon).$$
	When the above two values coincide,
	it is called the \textbf{Minkowski dimension} of $(X,d)$ and $\mathrm{is~denoted~by~dim}_{\mathrm{M}}(X,d).$  
	
	Writing
	$$\text{dim} _ {\text{M}} (X,d_{F},\varepsilon)=\frac{\log\#(X,d_{F},\varepsilon)}{\log(1 /\varepsilon)}.$$
	
	For every $\varepsilon>0$, we define 
	$$S(X,G,d,\varepsilon)=\lim_{n \rightarrow \infty}\frac{\log\#(X,d_{F_{n}},\varepsilon)}{|F_{n}|}.$$
	
	It's clear that $h(F):=\log\#(X,d_{F},\varepsilon)$ satisfies Lemma \ref{l100} and hence the limit always exists and does not depend on the choice of the F$\phi$lner sequence $\{F_{n}\}^{\infty}_{n=1}$.
	The \textbf{topological entropy} of $(X,G)$ is defined by
	$$h_{\mathrm{top}}(X,G)=\lim_{\varepsilon \rightarrow 0}S(X,G,d,\varepsilon).
	$$
	
	The \textbf{upper and lower metric mean dimensions} are defined by
	$$
	\overline{\text{mdim}} _ {\text{M}} (X,G,d)=\limsup_{\varepsilon \rightarrow 0 }\frac{S(X,G,d,\varepsilon)}{\log(1 /\varepsilon)}=\limsup_{\varepsilon \rightarrow 0 }\left(\lim_{n \rightarrow \infty}\frac{\text{dim} _ {\text{M}} (X,d_{F_{n}},\varepsilon)}{|F_{n}|}\right),
	$$
	$$
	\underline{\text{mdim}} _ {\text{M}} (X,G,d)=\liminf_{\varepsilon \rightarrow 0 }\frac{S(X,G,d,\varepsilon)}{\log(1 /\varepsilon)}=\liminf_{\varepsilon \rightarrow 0 }\left(\lim_{n \rightarrow \infty}\frac{\text{dim} _ {\text{M}} (X,d_{F_{n}},\varepsilon)}{|F_{n}|}\right).
	$$
	When the above two values coincide,
	it is called the \textbf{metric mean  dimension} of $(X,G)$ and $\mathrm{is~denoted~by~mdim}_{\mathrm{M}}(X,G,d).$  
	\subsection{Mean Hausdorff dimension}
	\
	\newline \indent \quad\ For $s\geq 0$ and $\varepsilon> 0$, we define $\mathcal{H} _\varepsilon^s( X, d)$ by
	$$\mathcal{H}_{\varepsilon}^{s}(X,d)=\inf\left\{\sum_{i=1}^{\infty}\left(\operatorname{diam}E_{i}\right)^{s}\Bigg|X=\bigcup_{i=1}^{\infty}E_{i}\text{ with diam }E_{i}<\varepsilon(\forall i\geq1)\Bigg\}.\right. $$
	
	\hspace{4mm}
	For $\varepsilon>0$, writing
	$$\dim_{\mathrm{H}}(X,d,\varepsilon)=\sup\{s\geq0\mid \mathcal{H}_\varepsilon^s(X,d)\geq1\}.$$
	The \textbf{Hausdorff dimension} of $(X,d)$ is defined by
	$$\dim_{\text{H}}(X,d)=\lim_{\varepsilon\to0}\dim_{\text{H}}(X,d,\varepsilon).$$
	For $0<\varepsilon<1,$ we have 
	$$\operatorname{dim}_{\mathrm{H}}(X,d,\varepsilon)\leq\dim_{\mathrm{M}}(X,d,\varepsilon).  $$
	Hence
	$$\dim_{\mathrm{H}}(X,d)\le\underline{\dim}_{\mathrm{M}}(X,d)\le\overline{\dim}_{\mathrm{M}}(X,d).$$

	\hspace{4mm}
	Let $\{F_{n}\}^{\infty}_{n=1}$ be a F$\phi$lner sequence in $G$. The \textbf{upper and lower mean Hausdorff dimensions} of $(X,G)$ with respect to $\{F_{n}\}^{\infty}_{n=1}$ are defined by
	\begin{align*}	\overline{\operatorname{mdim}}_{\mathrm{H}}(X,\{F_{n}\},d)&=\lim_{\varepsilon\to0}\left(\operatorname*{limsup}_{n\to\infty}\frac{\operatorname{dim}_{\mathrm{H}}\left(X,d_{F_{n}},\varepsilon\right)}{|F_{n}|}\right),\\
		\underline{\operatorname{mdim}}_{\mathrm{H}}(X,\{F_{n}\},d)&=\lim_{\varepsilon\to0}\left(\operatorname*{liminf}_{n\to\infty}\frac{\operatorname{dim}_{\mathrm{H}}\left(X,d_{F_{n}},\varepsilon\right)}{|F_{n}|}\right).
	\end{align*}
	These values depend on the choice of the F$\phi$lner sequence $\{F_{n}\}^{\infty}_{n=1}$. 
	When the above two values coincide, the common value is called the \textbf{mean Hausdorff dimension} with respect to $\{F_{n}\}^{\infty}_{n=1}$ and denoted by $\mathrm{mdim}_{\mathrm{H}}(X,\{F_{n}\},d).$
	
	\hspace{4mm}
	By  $||\cdot||_{\infty}$ we denote the $\ell^{\infty}$-norm on $\R^{N}$. To get the comparative relationship  between the mean dimension and mean Hausdorff dimension of $G$-action, we need some preliminary lemmas.
	\begin{lem}\cite{LT1} \label{X300}
		There exist a positive number $L$, a positive integer $M$ and an $L$-Lipschitz map $\varphi: (K,d) \rightarrow ([0,1]^{M}, ||\cdot||_{\infty})$ such that if $x,y \in K$ satisfy
		$||\varphi(x)-\varphi(y)||_{\infty} <1$, then $d(x,y) < \varepsilon$.
	\end{lem}

	\begin{lem}\cite{LT1}\label{X100}
		Let $N$ be a positive integer and $\varepsilon, \delta, s, \tau, L$ positive numbers with $4^{N}L^{s+\tau}\delta^{\tau} <1$. Let $(X,d)$ be a compact metric space with $\mathrm{dim}_{\mathrm{H}}(K, d, \delta) <s$. Suppose there exists an $L$-Lipschitz map $\varphi: (K,d) \rightarrow ([0,1]^{N}, ||\cdot||_{\infty})$ such that if $x,y \in K$ satisfy $||\varphi(x)-\varphi(y)||_{\infty} <1$ then $d(x,y) < \varepsilon$. Then $\mathrm{Widim}_{\varepsilon}(K,d) \leq s+ \tau$.
	\end{lem}
	
	\hspace{4mm}
	The following theorem shows the relationship of mean dimension, mean Hausdorff dimension and metric mean dimension for the amenable group actions.
	
	\begin{thm}\label{x000}
		Let $G$ be a countable discrete amenable group and $(X, G)$ be a $G$-system with a metric $d$. If $\{F_{n}\}_{n=1}^{\infty}$ is a F$\phi$lner sequence in $G$, then
		$$\mathrm{mdim}(X, G) \leq \underline{\mathrm{mdim}}_{\mathrm{H}}(X, \{F_{n}\}, d)\le\overline{\mathrm{mdim}}_{\mathrm{H}}(X,\{F_{n}\},d)      
		\le\underline{\mathrm{mdim}}_{\mathrm{M}}(X,G,d).
		$$
	\end{thm}
	
	\begin{proof}
		We only show $\mathrm{mdim}(X, G) \leq \underline{\text{mdim}}_{\text{H}}(X, \{F_{n}\}, d)$, the remaining inequalities are obvious.
		Let $\varepsilon>0$. We  can assume $\underline{\mathrm{mdim}}_{\text{H}}(X,\{F_{n}\},d) < \infty$. We take $\tau >0$ and $s>\underline{\mathrm{mdim}}_{\text{H}}(X,\{F_{n}\},d)$. Let $n \geq 1$ and define an $L$-$Lipschitz$  map $$\varphi_{F_{n}}: (X,d_{F_{n}}) \rightarrow \left(([0,1]^{M})^{F_{n}},||\cdot||_{\infty}\right)$$
		by $\varphi_{F_{n}}(x)=(\varphi(gx))_{g \in F_{n}}$, where $\varphi$ is defined in Lemma \ref{X300} and $\{F_{n}\}$ is a F$\phi$lner sequence. By Lemma \ref{X300}, if $x,y \in X$ satisfy $||\varphi_{F_{n}}(x)-\varphi_{F_{n}}(y)||_{\infty} <1$, then $d_{F_{n}}(x,y) < \varepsilon$.
		
		\hspace{4mm}
		Choose a  $\delta>0$ satisfying $4^{M}L^{s+\tau}\delta^{\tau} <1$. It follows from $\mathrm{mdim}_{\text{H}}(X,\{F_{n}\},d)<s$ that there exist $n_{1} < n_{2}<n_{3}< \cdots \rightarrow \infty$ satisfying $\mathrm{dim}_{\text{H}}(X,d_{F_{n}},\delta)<s|F_{n_{k}}|$. Note that $$4^{M|F_{n_{k}}|}L^{s|F_{n_{k}}|+\tau |F_{n_{k}}|}\delta^{\tau |F_{n_{k}}|}=(4^{M}L^{s+\tau}\delta^{\tau})^{|F_{n_{k}}|}<1.$$
		Then we can apply Lemma \ref{X100} to the space $(X,d_{F_{n_{k}}})$ and the map $\varphi_{F_{n_{k}}}$, which implies that 
		$$\mathrm{Widim}_{\varepsilon}(X,d_{F_{n_{k}}}) \leq s|F_{n_{k}}|+ \tau|F_{n_{k}}|.
		$$
		Therefore 
		$$\lim_{n \rightarrow \infty}\frac{\mathrm{Widim}_{\varepsilon}(X,d_{F_{n}})}{|F_{n}|} \leq s+ \tau.
		$$
		The right-hand side is independent of $\varepsilon$. Thus $\mathrm{mdim}(X,G) \leq s+\tau$. Let $\tau \rightarrow 0$ and $s \rightarrow \underline{\mathrm{mdim}}(X,\{F_{n}\},d)$, we obtain desired statement.
	\end{proof}
	
	\hspace{4mm}
	To better illustrate the relationship between metric mean dimension and mean Hausdorff dimension, we will give examples of their equivalence and inequality, respectively, in the Appendix.
	
	\subsection{Weighted topological entropy}
	\
	\newline \indent \quad\ Weighted topological entropy was introduced by Feng and Huang \cite{FH}, followed by Tsukamoto, who provided an equivalent definition in \cite{TM1}. Subsequently, Yang et al. studied this concept in the context of amenable group actions \cite{YCY}. Suppose $(X, G, T)$  and  $(Y, G, S)$  are two $G$-systems, and $d$  and  $d^{\prime}$  are metrics on  $X$ and $Y$, respectively. If there is a continuous surjection $\pi: X \rightarrow Y$ such that $\pi (gx)=g (\pi x)$ for any $g\in G, x\in X$, then we call $\pi$ a factor map between the $G$-systems $(X,G,T)$ and $(Y,G.S)$.
	
	\hspace{4mm}
	Let  $\pi: X \rightarrow Y$ be a factor map and $0 \leq w \leq 1$. Let  $d$  and  $d^{\prime}$  be metrics on  $X$  and  $Y$,   respectively. For each natural number $n$, we define new metrics  $d_{F_{n}}$  and  $d_{F_{n}}^{\prime}$  on  $X$  and  $Y$ by
	$$d_{F_{n}}\left(x_{1}, x_{2}\right)=\max _{g\in F_{n}} d\left(T^{g} x_{1}, T^{g} x_{2}\right), d_{F_{n}}^{\prime}\left(y_{1}, y_{2}\right)=\max _{g\in F_{n}} d^{\prime}\left(S^{g} y_{1}, S^{g} y_{2}\right).$$
	 Writing $\#^{w}\left(\pi, X, d_{F_{n}}, d_{F_{n}}^{\prime}, \varepsilon\right)$  
	$$\begin{array}{l}
		\#^{w}\left(\pi, X, d_{F_{n}}, d_{F_{n}}^{\prime}, \varepsilon\right) \\
		=\min \left\{\sum_{k=1}^{t}\left(\#\left(\pi^{-1}\left(V_{k}\right), d_{F_{n}}, \varepsilon\right)\right)^{w} \mid \begin{array}{c}
			\{V_{1} \cup \cdots \cup V_{n}\} \text { is an open cover of } Y \text{ with }\\
			\operatorname{diam}\left(V_{k}, d_{F_{n}}^{\prime}\right)<\varepsilon \text { for all } 1 \leq k \leq n
		\end{array}\right\} . \\
	\end{array}$$
	
	\hspace{4mm}
	It is easy to check that $\Omega \in Fin(G) \rightarrow \log\#^{w}\left(\pi, X, d_{\Omega}, d_{\Omega}^{\prime}, \varepsilon\right)$ satisfies the conditions in Lemma \ref{l100}. Then
	by Lemma \ref{l100}, the limit 
	$\lim _{n \rightarrow \infty} \frac{\log \#^{w}\left(\pi, X, d_{F_{n}}, d_{F_{n}}^{\prime}, \varepsilon\right)}{|{F_{n}}|}$
	exists. Then the  \textbf{$w$-weighted topological entropy} of the factor map  $\pi: X \rightarrow Y$ is defined by
	$$h_{\text {top }}^{w}(\pi, G)=\lim _{\varepsilon \rightarrow 0}\left(\lim _{n \rightarrow \infty} \frac{\log \#^{w}\left(\pi, X, d_{F_{n}}, d_{F_{n}}^{\prime}, \varepsilon\right)}{|{F_{n}}|}\right).$$ 
	
	\subsection{Topological Entropy for \texorpdfstring{$G \times \N$}{}-actions}
	\
	\newline \indent \quad\ If $X$ is a compact metric space, $G$ is a  countable discrete amenable group, and $H:G\times X\rightarrow X, (g,x)\mapsto H^gx$ and $T:\N\times X\rightarrow X,(n,x)\mapsto T^nx$   are continuous maps with  $H^g T^{n}(x)=T^{n}H^g(x)$ for any $g\in G$ and any $n \in \mathbb{N}$, then we call  $(X,G,H,T)$ a $G \times \N$-actions.

	\hspace{4mm}
	Let $ (X, G, H, T)$  be a  $G \times \N$-actions with a metric  $d$  on  $X$. Given a subset  $\Omega \subset G \times \N$, we can define a metric  $d_{\Omega}^{H, T}$  on $ X$  by
	$$d_{\Omega}^{H, T}(x, y)=\sup _{(g, w) \in \Omega} d\left(H^{g} T^{w}( x), H^{g} T^{w} (y)\right).$$ 
	
	\hspace{4mm}
	Let $\{F_{n}\}_{n=1}^{\infty}$ be a F$\phi$lner sequence in $G$ and $F^{'}_{n}=\{0,1,...,n-1\} \subset \N$. We define 
	\begin{align*}
		d_{F_{n} \times F^{'}_{n}}^{H, T}(x, y) & =\max _{\substack{g\in F_{n} ,
				w\in F^{'}_{n}}} d\left(H^{g} T^{w}( x), H^{g} T^{w} (y)\right),  
	\end{align*}
	
	and define the map $h:\Omega \in Fin(G \times \N)  \rightarrow \mathbb{R}$ by 
	$$h(\Omega)=\log\#\left(X, d_{\Omega}^{H, T}, \varepsilon\right).$$
	Obviously, the map $h$ satisfies the conditions in Lemma \ref{l100}.
	Let $F_{n} \times F^{\prime}_{n}=\{(u,v)\mid u\in F_{n}\text{ and } v\in F^{\prime}_{n}\}$, then $\{F_{n} \times F^{\prime}_{n}\}_{n \geq 1}$ is a F$\phi$lner sequence in $G \times \N$ and $|F_{n}\times F^{\prime}_{n}|=n\cdot |F_{n}|$.
	By Lemma \ref{l100}, the limit
	$\lim _{n \rightarrow \infty} \frac{\log\#\left(X, d_{F_{n} \times F^{\prime}_{n}}^{H, T}, \varepsilon\right)}{n\cdot|{F_{n}}|}$
	exists. We define the \textbf{topological entropy} of  $(X, G \times \N, H, T)$ by
	$$h_{\text {top }}(X, G \times \N, H, T) =\lim _{\varepsilon \rightarrow 0}\left(\lim _{n \rightarrow \infty} \frac{\log \#\left(X, d_{F_{n} \times F^{\prime}_{n}}^{H, T}, \varepsilon\right)}{n\cdot | F_{n}|}\right).$$
	We also have
	$$h_{\text {top }}(X, G \times \N, H, T)=\lim _{\varepsilon \rightarrow 0}\left(\lim _{\substack{m \rightarrow \infty \\ n \rightarrow \infty}} \frac{\log \#\left(X, d_{F_{n}\times F^{\prime}_{m}}^{H, T}, \varepsilon\right)}{m\cdot |F_{n}|}\right).$$
	\section{Self-similar systems and Homogeneous systems for amenable group actions}
	\hspace{4mm}
	In this section, we will study the relationship of mean Hausdorff dimension and metric mean dimension of self-similar systems and homogeneous systems, we find that the mean Hausdorff dimension and metric mean dimension of these two systems are equal. We will start with introducing self-similar systems for amenable group actions.
	
	\subsection{Self-similar Systems}	
	\
	\newline \indent \quad\ Let $G$ be a countable discrete amenable group and $\|x\|_\infty$ be the $\ell^\infty$-norm on $\mathbb{R}^{G}$. Assuming that $\ell^{\infty}$ is endowed with the weak$^*$ topology. We define the shift map $\sigma:G \times \ell^\infty\to\ell^\infty$ by
	$$
	\sigma^{h} \left((x_{g})_{g\in G}\right)=(x_{gh})_{g \in G},\text{ for any } h \in G .
	$$
	This is continuous concerning the weak$^*$ topology. We will consider a certain self-similar set of $\ell^{\infty}$ which is $\sigma$-invariant.
	
	\hspace{4mm}
	Let $( \Omega,G, T) $ be a $G$-system, where  the continuous map  $T:G \times \Omega \to \Omega$ is defined by $(g,x) \mapsto  gx$. Suppose that for each $\omega \in\Omega$, there is a point $H( \omega) = \left ( H( \omega) _{g}\right ) _{g\in G }\in \ell^{\infty}$ such that the map $H:\Omega \rightarrow \ell^{\infty}, \omega \mapsto H(\omega)$
	is continuous with respect to the weak* topology of $\ell^\infty$, and $\sigma^{g} (H(\omega))=$ $H( T^{g} (\omega)), \text{ for any } g \in G \text{ and any } \omega \in \Omega.$ Since $\Omega$ is compact, we have
	$$
	\sup_{\omega\in\Omega}\left\|H(\omega)\right\|_\infty<\infty.
	$$
	Fix $c \in \R$ with $0<c<1.$ For every $\omega\in\Omega$, we can define a contracting similarity transformation $S_{\omega}: \ell^{\infty}\to \ell^{\infty}$ by
	$$
	S_\omega(x)=cx+H(\omega).
	$$
	Then for any $g \in G$, we have
	$$\sigma^{g}\left ( S_\omega( x) \right ) = S_{T^{g}\omega}\left ( \sigma^{g}( x) \right ) .$$  
	
	\hspace{4mm}
	Now we will define the self-similar system on amenable group actions.
	\begin{definition}
		Let $Y$ be a nonempty subset of $\ell^\infty$ and $F:G \times \ell^\infty\to\ell^\infty$ be a shift map. We call  $(Y, G, F)$ a \textbf{self-similar system} if $Y=\cup_{i \in \Omega}S_i(Y)$, where $\{S_i\}_{i \in \Omega}$ is a family of contracting similarity transformations and $F^{g}(Y)\subset Y$ for any $g \in G$.
	\end{definition}
	
	\begin{pro} There uniquely exists a self-similar dynamical system $(X,G,\sigma)$ in $\ell^{\infty}$, where $$
		X=\bigcup_{\omega\in\Omega}S_\omega(X),
		$$
		and $\{{S_\omega}\}_{\omega\in\Omega}$ is a family of contracting similarity transformations. 
	\end{pro}
	
	\begin{proof}
		For the proof one may refer to Proposition 4.1 of \cite{TM}, we omit it here.
	\end{proof}
	
	\subsection{The calculation for self-similar systems}
	\
	\newline \indent \quad\ Now we are going to prove the relationship of mean Hausdorff dimension and metric mean dimension of the above self-similar system, we will begin with some preparations. The following result from Tsukamoto \cite{LT2} serves as an important lemma for computation and is a finite accuracy version of a theorem by Falconer \cite{FKJ} 
	
	\begin{lem}\cite{TM}\label{l46}
		Given $0<\varepsilon,a,\tau<1$ and a compact metric space (X,d). Suppose for any closed ball $B\subset X$ of radius $\varepsilon$, there exists a map $\varphi:X\to B$ satisfying
		$$
		d\left(\varphi(x),\varphi(y)\right)\geq a\varepsilon d(x,y)\quad(x,y\in X).
		$$
		Then we have
		$$\dim_{\mathrm{H}}(X,d,\delta_0)\geq\tau\cdot\frac{\log \#(X,d,9\varepsilon)}{\log\left(\frac{1}{a\varepsilon}\right)}.$$
	\end{lem}
	
	\hspace{4mm}
	We recall that $X$ is compact. Therefore, we can find $D\geq1$ such that
	$$
	\mathrm{diam}(X,d)<D.
	$$
	For $n \geq 1$, we define a metric $d_{F_{n}}$ on $X$ by
	$$
	d_{F_{n}}(x,y)=\max_{h \in F_{n}}d\left(\sigma^h(x),\sigma^h(y)\right).
	$$
	\begin{claim} \label{c47}
		Let $n$ be a natural number and let $0<\varepsilon<1$. For any point $ q \in X $, there exists a map $\varphi:X\to \overline{B}_{d_{F_{n}}}(q,\v)$ that satisfys
		$$		d_{F_{n}}\left(\varphi(x),\varphi(y)\right)\geq\frac{c\varepsilon}Dd_{F_{n}}(x,y).
		$$
		Here $\overline{B}_{d_{F_{n}}}(q,\varepsilon)$ is the closed $\varepsilon$-ball centered at $q$ with respect to the metric $d_{F_{n}}.$ 
	\end{claim}
	\begin{proof}
		Given a point $\omega\in\Omega$, we have 	$$d\left(S_\omega(x),S_\omega(y)\right)=d\left(cx+H(\omega),cy+H(\omega)\right)=\sum_{g \in G}\alpha_{g}|cx_{g}-cy_{g}|=cd(x,y). $$
		Similarly, we have $ d_{F_{n}}\left(S_\omega(x),S_\omega(y)\right)=cd_{F_{n}}(x,y). $
		Choose $k\in \N$ such that $c^kD\leq\varepsilon<c^{k-1}D$. Since $X=\bigcup_{\omega\in\Omega}S_\omega(X)$, we can find a sequence $\omega_1,\ldots,\omega_k\in\Omega$  satisfying
		$$
		q\in S_{\omega_1}\circ\cdots\circ S_{\omega_k}(X).
		$$
		We define a set $\varphi:=S_{\omega_1}\circ\cdots\circ S_{\omega_k}:X\to X.$ Then for any $x,y\in X,$ we have
		$$
		d_{F_{n}}(\varphi(x),\varphi(y))=c^kd_{F_{n}}(x,y)\leq c^kD\leq\varepsilon,
		$$
		Hence $\varphi(X)\subset B_{d_{F_{n}}}(q,\varepsilon).
		$ Since $\varepsilon<c^{k-1}D$, 
		it follows that $c^k>\frac{c\varepsilon}{D}$.  Therefore,\begin{equation*}
			d_{F_{n}}\left(\varphi(x),\varphi(y)\right)=c^kd_{F_{n}}(x,y)\ge\frac{c\varepsilon}{D}d_{F_{n}}(x,y).\qedhere
		\end{equation*}
	\end{proof}
	
	\hspace{4mm}
	The following conclusion shows the relationship between metric mean dimension and mean Hausdorff dimension of self-similar systems.
	
	\begin{thm}\label{t44}
		Let $(X,G,\sigma)$ and $(\Omega,G,T)$ be the above two $G$-systems. Suppose $\{F_{n}\}$ is a F$\phi$lner sequence in $G$, then  
		$$  \mathrm{mdim}_{\mathrm{M}}(X,G,d)=\mathrm{mdim}_{\mathrm{H}}(X,\{F_{n}\},d).$$
	\end{thm}
	\begin{proof}  By Lemma \ref{x000}, we have 
		$$\underline{\mathrm{mdim}}_{\mathrm {H}}(X,\{F_{n}\},d)\leq\overline{\mathrm {mdim}}_{\mathrm {M}}(X,G,d).$$
		We will prove the reverse inequality.
		
		\hspace{4mm}
		Let $0<\varepsilon,\tau<1.$ By Lemma \ref{l46} and Claim \ref{c47}, we can find $\delta_0=\delta_0(\varepsilon,c/D,\tau)>0$ such that for any $n\in \N$,
		$$\dim_{\mathrm{H}}(X,d_{F_{n}},\delta_0)\geq\tau\cdot\frac{\operatorname{log}\#(X,d_{F_{n}},9\varepsilon)}{\operatorname{log}\left(\frac{D}{c\varepsilon}\right)}.$$
		Then
		$$
		\underline{\mathrm{mdim}}_{\mathrm{H}}(X,\{F_{n}\},d)\geq\liminf_{n\to\infty}\frac{\dim_{\mathrm{H}}(X,d_{F_{n}},\delta_0)}{|F_{n}|}\geq \lim_{n\to\infty}\tau\cdot\frac{\log\#(X,d_{F_{n}},9\varepsilon)}{|F_{n}|\log\left(\frac D{c\varepsilon}\right)}.
		$$
		Letting $\varepsilon \rightarrow 0$ and $\tau \rightarrow1$, we get
		\begin{equation*}	    
			\underline{\mathrm{mdim}}_{\mathrm {H}}(X,\{F_{n}\},d)\geq\overline{\mathrm {mdim}}_{\mathrm {M}}(X,G,d).\qedhere
		\end{equation*}
	\end{proof}

 \hspace{4mm}
 Let $E$ be a compact subset of $X$, $F \in Fin(G)$ and $\varepsilon>0$. Recall that $K \sub  E $ is an $(F,\varepsilon)$-spanning set of $E$ if  for any $x \in E$ there exists $y \in K$ such that $d_{F}(x,y) \leq \varepsilon$.
 
	\hspace{4mm}
	The following Theorem shows the relationship between the topological entropy and metric mean dimension of self-similar systems.
	
	\begin{thm}\label{p3}
		For the above two $G$-systems $(X,G,\sigma)$, we have 
		$$ \mathrm{mdim}_{\mathrm{M}}(X,G,d)\leq \frac{h_{\mathrm{top}}(\Omega,G)}{\log(1/c)}.$$
		Here $h_{\mathrm{top}}(\Omega,G)$ is the topological entropy of $(\Omega,G,T)$.
	\end{thm}
	\begin{proof}
		Assume that $ \rho$ is a metric on $\Omega.$ For $n\geq1$ we define a metric $\rho_{F_{n}}$ on $\Omega$ by
		$$\rho_{F_{n}}(\omega,\omega')=\max_{g \in F_{n}}\rho\left(T^g\omega,T^g\omega'\right), \text{for any } \omega,\omega' \in \Omega.$$
		Given $\omega,\omega'\in\Omega$. For any $x,y\in X$, we have 
		\begin{align*}
			d\left(S_{\omega}(x),S_{\omega'}(y)\right)
			&=\sum_{g \in  G}\alpha_{g}\left|cx_g+H(\omega)_g-cy_g-H(\omega^{\prime})_g\right| \\
			&\leq cd(x,y)+d\left(H(\omega),H(\omega')\right).
		\end{align*}
		
		Then for any $n\in \N$,
		$$d_{F_{n}}\left(S_\omega(x),S_{\omega^{\prime}}(y)\right) \leq cd_{F_{n}}(x,y)+d_{F_{n}}\left(H(\omega),H(\omega^{\prime})\right).$$
		Recall that $\mathrm{diam}( X, d_{F_{n}}) \leq \mathrm{diam}( X, d) < D$. By repeatedly applying the inequality stated above for ${\omega_1,\cdots,\omega_n,\omega_1',\cdots,\omega_n'\in\Omega }$ and $ {x,y\in X}$, we obtain
		\begin{align*}
			&d_{F_{n}}\left(S_{\omega_1}\circ\cdots\circ S_{\omega_n}(x),S_{\omega_1^{\prime}}\circ\cdots\circ S_{\omega_n^{\prime}}(y)\right)\\
			&\leq c^nd_{F_{n}}(x,y)+\sum_{i=1}^nc^{i-1}d_{F_{n}}\left(H(\omega_i),H(\omega_i')\right)\\
			&\leq c^nd_{F_{n}}(x,y)+\frac{\max_{1\leq i\leq n}d_{F_{n}}\left(H(\omega_i),H(\omega_i^{\prime})\right)}{1-c} \\
			&<c^nD+\frac{\max_{1\leq i\leq n}d_{F_{n}}\left(H(\omega_i),H(\omega_i^{\prime})\right)}{1-c}.
		\end{align*}

		\hspace{4mm}
		Recall that $\Omega$ is a compact metric space. For any $0< \varepsilon< 1,$ we can choose $\delta>0$ such that if $\rho(\omega,\omega^{\prime})<\delta$, then $d\left(H(\omega),H(\omega^{\prime})\right)<(1-c)\frac\varepsilon6$. Similarly, for any $n \in \N$
		if $
		\rho_{F_{n}}(\omega,\omega^{\prime})=\max_{g \in F_{n}}(T^{g}(\omega),T^{g}(\omega^{\prime}))<\delta$, then 
		\begin{align*}
			d_{F_{n}}\left(H(\omega),H(\omega^{\prime})\right)&=\max_{g \in F_{n}}d(\sigma^{g}(H(\omega)),\sigma^{g}(H(\omega^{\prime}))\\
			&=\max_{g \in F_{n}}d(H(T^{g}(\omega)),H(T^{g}(\omega^{\prime})))\\
			&<(1-c)\frac\varepsilon6.
		\end{align*}
		
		\hspace{4mm}
		We also take $m\in \N$  such that
		$$Dc^m <\frac{\varepsilon}{6}\leq Dc^{m-1}.$$
		
		\hspace{4mm}
		Choose a point $p\in X$. For $n \geq m$, let $\Omega_{n,\delta}\subset\Omega$ be a $\delta$-spanning set with respect to $\rho_{F_{n}}$ with $|\Omega_{n,\delta}|=\#(\Omega,\rho_{F_{n}},\delta)$. For any $\omega_1,\cdots,\omega_m\in\Omega$, there exists $\omega_1',\cdots,\omega_m'\in\Omega_{n,\delta}$ satisfying $\rho_{F_{n}}(\omega_i,\omega_i')<\delta$, then $d_{F_{n}}(H(\omega_i),H(\omega_i')<(1-c)\varepsilon/6.$ So for any $x\in X,$ we have
		$$d_{F_{n}}(S_{\omega_1}\circ\cdots\circ S_{\omega_m}(x),S_{\omega_1'}\circ\cdots\circ S_{\omega_m'}(p))<\frac{\varepsilon}{3}.$$
		
		Since	$X=\bigcup_{\omega_1,\cdots,\omega_m\in\Omega}S_{\omega_1}\circ\cdots\circ S_{\omega_m}(X)$, this implies
		$$\left\{\left.S_{\omega_1^{\prime}}\right.\circ\cdots\circ S_{\omega_m^{\prime}}(p)\right|\left.\omega_1^{\prime},\cdots,\omega_m^{\prime}\in\Omega_{n,\delta}\right\}$$
		is a $\varepsilon/3$-spanning set of $X$ with respect to $d_{F_{n}}$. Then we have $$\#\left(X,d_{F_{n}},\varepsilon\right)\leq |\Omega_{n,\delta}|^m=\{\# (\Omega,\rho_{F_{n}},\delta)\}^m.$$
		Hence
		$$\frac{\log\#\left(X,d_{F_{n}},\varepsilon\right)}{|F_{n}|\log\left(\frac{6D}{c\varepsilon}\right)}\le\frac{\log\#\left(\Omega,\rho_{F_{n}},\delta\right)}{|F_{n}|\log(1/c)}, $$
		which implies that
		\begin{align*}
			\mathrm{mdim}_{\mathrm{M}}(X,G,d)
			&=\lim_{\varepsilon \rightarrow 0 }\left(\lim_{n \rightarrow \infty}\frac{\text{dim} _ {\text{M}} (X,d_{F_{n}},\varepsilon)}{|F_{n}|}\right)\\
			& \leq \lim_{\varepsilon \rightarrow 0 }\left(\lim_{n\to\infty}\frac{\log\#\left(X,d_{F_{n}},\varepsilon\right)}{|F_{n}|\log\left(\frac{6D}{c\varepsilon}\right)}\right)\\
			& \leq\lim_{\varepsilon \rightarrow 0 }\left( \lim_{n\to\infty}\frac{\log\# \left (\Omega,\rho_{F_{n}},\delta\right)}{|F_{n}|\log(1/c)}\right)\\
			&\leq\frac{h_{\mathrm{top}}(\Omega,G)}{\log(1/c)}.\qedhere
		\end{align*}
	\end{proof}
	
	\hspace{4mm}
	Combining Theorem \ref{t44} and Theorem \ref{p3}, we obtain Theorem \ref{t1.2}.

	\subsection{Homogeneous systems} 
	\
	\newline \indent \quad\ Let $G$ be a countable discrete amenable group. Suppose $\mathbb{R} / \mathbb{Z}$  is a circle with a metric  $\rho$  defined by
	$$\rho(x, y)=\min _{n \in \mathbb{Z}}|x-y-n|,$$
	and $(\mathbb{R} / \mathbb{Z})^{G}$ is the infinite-dimensional torus.
	A metric  $d$ of $(\mathbb{R} / \mathbb{Z})^{G}$ is given by
	$$d\left(\left(x_{g}\right)_{g \in G},\left(y_{g}\right)_{g \in G}\right)=\sum_{g\in G} \alpha_{g} \rho\left(x_{g}, y_{g}\right) ,$$
	where $\alpha_{g} \in (0,+\infty)$ satisfies $$
	\alpha_{1_{G}}=1,\sum_{g \in G}\alpha_{g} < +\infty.
	$$
	
	The shift map $ \sigma: G\times (\mathbb{R} / \mathbb{Z})^{G} \rightarrow(\mathbb{R} / \mathbb{Z})^{G} $ is defined by 
	$$\sigma^h\left(\left(x_{g}\right)_{g \in G}\right)=\left(x_{gh}\right)_{g \in G}, \text{ for any } h\in G.$$
	Let  $ b\in \N$  with $b>1$. The map $T_b: \N\times (\mathbb{R} / \mathbb{Z})^{G} \rightarrow(\mathbb{R} / \mathbb{Z})^{G}$ is defined by $$T_b^m((x_g)_{g\in G})= (b^mx_g)_{g\in G}.$$
	
	\hspace{4mm}
	If there is a closed subset $X\subset(\mathbb{R} / \mathbb{Z})^{G}$ satisfies that $T^{m}_b(X) \subset X$ for any $m\in \N$ and $ \sigma^{h}(X) \subset X$ for any $h\in G,$ then we call $X$ a \textbf{homogeneous system}.
	
	\subsection{The calculation for the homogeneous systems} 
        \
	\newline \indent \quad\ We will now begin calculating the mean Hausdorff dimension and metric mean dimension for the aforementioned homogeneous systems. The topological entropy for \texorpdfstring{$G \times \N$}{}-actions will play a crucial role in this process. First, we need to establish the following auxiliary tool to provide a lower bound for the mean Hausdorff dimension of homogeneous systems
	\begin{lem}\cite{TM}\label{l1}
		Suppose  $X$  is a compact metric space with a metric  $d $,  and $T : X \rightarrow X$ is a Lipschitz continuous map with Lipschitz constant $L >1$. If there are positive numbers $t, \delta, \varepsilon$ that satisfy
		$$0<t<1, \quad 0<\delta<1, \quad \delta^{1-t}<\varepsilon,$$
		then 
		$$\lim _{N \rightarrow \infty} \frac{\log \#\left(X, d_{N}^{T}, \varepsilon\right)}{N} \leq \frac{\log L}{t} \cdot \operatorname{dim}_{\mathrm{H}}(X, d, \delta) .$$
	\end{lem}
	
	\begin{thm}\label{t11}
		If $(X,G, S, T) $ is a $ G\times \N$-actions on the space $X$ equipped with a metric $d$, and there exists $L>1 $ satisfying 	
		$$d(Tx, Ty) \leq L \cdot d(x, y) \text{ for any } x, y \in X,$$
		then we have
		$$\underline{\operatorname{mdim}}_{\mathrm{H}}(X, \{F_{n}\}, d) \geq \frac{h_{\mathrm {top }}(X, G \times \N, S, T)}{\log L} .$$
	\end{thm}
	
	\begin{proof}
		It is obvious that $ d_{F_{n}}^{S}(Tx, Ty) \leq L \cdot d_{F_{n}}^{S}(x, y) $ for  any natural number $n$. Let $0<\delta,t,\varepsilon<1$ be positive numbers with $\delta^{1-t}<\varepsilon.$ 
		
		\hspace{4mm}
		From Lemma \ref{l1}, for any natural number $ n>0 $, one has		
		$$\lim _{N \rightarrow \infty} \frac{\log \#\left(X, d_{F_{n}\times F^{'}_{N}}^{S, T}, \varepsilon\right)}{N} \leq \frac{\log L}{t} \cdot \operatorname{dim}_{\mathrm{H}}\left(X, d_{F_{n}}^{S}, \delta\right) .$$

		Then
		$$\lim _{\substack{ N \rightarrow \infty}} \frac{\log \#\left(X, d_{F_{n}\times F^{'}_{N}}^{S, T}, \varepsilon\right)}{N\cdot |F_{n}|} \leq \frac{\log L}{t} \cdot \frac{\operatorname{dim}_{\mathrm{H}}\left(X, d_{F_{n}}^{S}, \delta\right)}{|F_{n}|} .$$	
		
		Letting $n \rightarrow \infty$, we have
		$$\lim _{\substack{N \rightarrow \infty \\ n \rightarrow \infty}} \frac{\log \#\left(X, d_{F_{n}\times F^{'}_{N}}^{S, T}, \varepsilon\right)}{N\cdot |F_{n}|} \leq \frac{\log L}{t} \cdot \liminf _{n \rightarrow \infty} \frac{\operatorname{dim}_{\mathrm{H}}\left(X, d_{F_{n}}^{S}, \delta\right)}{|F_{n}|} ,$$		
		which implies	
		$$\lim _{\substack{N \rightarrow \infty \\ n \rightarrow \infty}} \frac{\log \#\left(X, d_{F_{n}\times F^{'}_{N}}^{S, T}, \varepsilon\right)}{N\cdot |F_{n}|} \leq \frac{\log L}{t} \cdot \underline{\operatorname{mdim}}_{\mathrm{H}}(X, \{F_{n}\}, d) .$$		
		Letting  $t \rightarrow 1$  and $ \varepsilon \rightarrow 0$ , we get the desire result \begin{equation*}\underline{\operatorname{mdim}}_{\mathrm{H}}(X, \{F_{n}\}, d) \geq \frac{h_{\text {top }}(X,G \times \N, S, T)}{\log L}. \qedhere
		\end{equation*}
	\end{proof}
	
	\hspace{4mm}
	In the following, we give an upper bound for the metric mean dimension of homogeneous systems. 
	\begin{pro}\label{p36}
		In the setting of Theorem \ref{t1.3}, we have
		$$\overline{\operatorname{mdim}}_{\mathrm{M}}(X, G,  d) \leq \frac{h_{\text {top }}(X, G \times \N, \sigma, T_b)}{\log b}.$$
	\end{pro}
	
	\begin{proof}
		Let $\varepsilon >0$. We choose a finite nonempty subset $S$ of $G$ such that 
		$$\sum_{g \in G\setminus S}\alpha_{g} \leq \frac{\varepsilon}{2}.$$
		Letting $c =\sum_{g \in G} \alpha_{g} < +\infty$. 
		
		\hspace{4mm}
		For any natural number $N>0$ and for any two points $ u, v \in \mathbb{R} / \mathbb{Z}$ satisfy $\max _{0\leq k< N} \rho\left(b^{k} u, b^{k} v\right)< \frac{1}{2cb}$, we have $\rho(u, v)< \frac{1}{2cb^{N}}$.
		Then for any $x, y\in  (\mathbb{R} / \mathbb{Z})^{G}$ with $x=(x_g)_{g\in G}, y=(y_g)_{g\in G}$, if we have  
		$$d_{S\times F^{'}_{N}}^{\sigma, T_b}(x, y)=\max_{(s,n) \in S \times F^{\prime}_{N}}d(\sigma^{s}T^{n}_{b}(x),\sigma^{s}T^{n}_{b}(y))<\frac{1}{2c b},$$
		which implies  for every  $g\in S$, we have	
		$$\max _{0\leq k<N} \rho\left(b^{k} x_{g}, b^k y_{g}\right)<\frac{1}{2c b},$$		
		then
		$$\rho\left(x_{g}, y_{g}\right)<\frac{1}{2cb^{N}}.$$
		We choose a natural number $N$ such that $b^{-N}\leq \varepsilon<b^{-N+1}$, then
		$$d(x, y)  \leq \sum_{g\in S} \alpha_{g} \rho\left(x_{g}, y_{g}\right)+\sum_{g \in G\setminus S} \alpha_{g}  <\frac{1}{2cb^{N}}\sum_{g\in S} \alpha_{g}+\sum_{g \in G\setminus S} \alpha_{g} <\frac{1}{2b^{N}}+\sum_{g \in G\setminus S} \alpha_{g}
		<\varepsilon.
		$$
		Hence for any natural number $n \geq 1$, $d_{SF_{n}\times F^{\prime}_{N}}^{\sigma, T_b}(x, y)<\frac{1}{2c b}$ implies that  
		$ d_{F_{n}}^{\sigma}(x, y)\leq \varepsilon$.	
		
		Therefore
		$$\#\left(X, d_{{F_{n}}}^{\sigma}, \varepsilon\right) \leq \#\left(X, d_{SF_{n}\times F^{\prime}_{N}}^{\sigma, T_b}, \frac{1}{2c b} \right).$$
		Then
		\begin{align*}
			\lim _{n \rightarrow \infty} \frac{\operatorname{dim}_{\mathrm{M}}\left(X, d_{F_{n}}^{\sigma}, \varepsilon\right)}{|F_{n}|} & =\lim _{n \rightarrow \infty} \frac{\log \#\left(X, d_{F_{n}}^{\sigma}, \varepsilon\right)}{|F_{n}| \log (1 / \varepsilon)}  \\
			&\leq \lim _{n \rightarrow \infty}\frac{\log \#\left(X, d_{SF_{n}\times F^{\prime}_N}^{\sigma, T_b}, \frac{1}{2c b}\right)}{|F_{n}| \log (1 / \varepsilon)}\\
			&= \lim _{n \rightarrow \infty}\frac{\log \#\left(X, d_{SF_{n}\times F^{'}_N}^{\sigma, T_b}, \frac{1}{2c b}\right)}{|SF_{n}| \log (1 / \varepsilon)} \cdot \frac{|SF_{n}|}{|F_{n}|}\\
			&= \lim _{n \rightarrow \infty}\frac{\log \#\left(X, d_{SF_{n}\times F^{'}_N}^{\sigma, T_b}, \frac{1}{2c b}\right)}{|SF_{n}| \log (1 / \varepsilon)} \cdot \lim _{n \rightarrow \infty}\frac{|SF_{n}|}{|F_{n}|}\\
			&\leq 	\lim _{n \rightarrow \infty}\frac{\log \#\left(X, d_{F_{n}\times F^{'}_N}^{\sigma, T_b}, \frac{1}{2c b}\right)}{ (N-1)|F_{n}| \log b}.\\
		\end{align*}
		Here we used a fact that $\lim_{n \rightarrow \infty}\frac{|SF_{n}|}{|F_{n}|}=1$ if $\{F_{n}\}_{n=1}^{\infty}$ is a F$\phi$lner sequence of $G$. Since $ N \rightarrow \infty$  as  $\varepsilon \rightarrow 0$, then 	
		\begin{align*}
			\limsup _{\varepsilon \rightarrow 0}\left(\lim _{n \rightarrow \infty} \frac{\operatorname{dim}_{\mathrm{M}}\left(X, d_{F_{n}}^{\sigma}, \varepsilon\right)}{|F_{n}|}\right) 
			&\leq\lim _{\substack{N \rightarrow \infty \\
					n \rightarrow \infty}}\frac{\log \#\left(X, d_{F_{n}\times F^{'}_N}^{\sigma, T_b}, \frac{1}{2c b}\right)}{|F_{n}|\cdot N \log b}\cdot \frac{N}{N-1}\\
			&\leq \frac{1}{\log  b} \cdot h_{\operatorname{top}}\left(X, G\times \N, \sigma, T_b\right).	
		\end{align*}
		This shows $\overline{\operatorname{mdim}}_{\mathrm{M}}(X,G, d) \leq \frac{h_{\operatorname{top}}\left(X, G \times \N, \sigma, T_b\right)}{\log b}$.
	\end{proof}
	
	\hspace{4mm}
	Combining  Theorem \ref{t11} and Proposition \ref{p36}, we obtain Theorem \ref{t1.3}.
	
	\section{Infinitely dimensional carpet of finitely generated amenable groups} 
	\hspace{4mm}
	In this section, we will introduce carpet systems of finitely generated amenable groups, we find that the metric mean dimension and mean Hausdorff dimension of infinitely dimensional carpets of finitely generated amenable groups are not equal. We will begin by reviewing some basic concepts related to finitely generated amenable groups, followed by an introduction to infinite-dimensional carpet systems for amenable group actions, and then proceed with our computations.
	
	\subsection{Finitely generated amenable groups}
	\
	\newline \indent \quad\ Let $G$ be an infinite discrete countable group. Recall that the subgroup generated by a subset $A \subset G$ is the smallest subgroup $ \langle A \rangle \subset G$ containing $A$. One says that group $G$ is \textbf{finitely generated} if there exists a finite subset $A \subset G$ such that $\langle A \rangle = G$. Let $G$ be a finitely generated amenable group with a symmetric generating set $S$. We say a generating set $S$ is symmetric if for any $s \in S$, we always have  $s^{-1} \in S$. The
	$S$-word-length $\ell_{S}(g)$ of an elements $g \in G$ is the minimal integer $n \geq 0$ such that $g$ can be expressed as a product of $n$ elements in $S$, that is ,
	$$
	\ell_{S}(g)=min\{n \geq 0: g=s_{1}\cdots s_{n}:s_{i} \in S, 1 \leq i \leq n\}.
	$$
	It immediately follows from the definition that $\ell_{S}(g)=0$ if and only if $g=1_{G}$. Define the metric $d_{S}$ on $G$: $d_{S}(g,h)=\ell_{S}(g^{-1}h)$. It is obvious that the metric $d_{S}$ is invariant by left multiplication.
	
	\hspace{4mm} 
	For $g \in G$ and $n \in \N$, we denote by
	$$
	B_{S}^{G}(g,n)=\{h \in G: d_{S}(g,h) \leq n\},
	$$
	the ball of radius $n$ in $G$ centered at the element $g \in G$. When $g=1_{G}$, we have $B_{S}^{G}(1_{G},n)=\{h \in G: \ell_{S}(h) \leq n\}$, and write $B_{S}^{G}(n)$ instead of $B_{S}^{G}(1_{G},n)$. Also, when there is no ambiguity on the group $G$, we omit the subscript $G$ and write $B_{S}(g,n)$ and $B_{S}(n)$ instead of $B_{S}^{G}(g,n)$ and $B_{S}^{G}(n)$. If $G$ is a finitely generated group with a symmetric generating set $S$, it's easy to know $\{B_{S}(m)\}_{m=1}^{\infty}$ is a F$\phi$lner sequence of $G$.
	
	\hspace{4mm}
	Let $G$ be a finitely generated group with a symmetric generating set $S=\{s_{1}, \cdots s_{m}\}$. We define an order in $S$ which formalizes through the following construction: given $s_{i}, s_{j} \in S$ we say that $s_{i}<s_{j}$ if $i <j$. Hence $s_{1}<s_{2}<\cdots<s_{m}$. This allows us to consider the order in $G$. For $g, g^{\prime} \in G,$ we call $g<g^{\prime}$ if $\ell_{S}(g)< \ell_{S}(g^{\prime})$. If $\ell_{S}(g)=\ell(g^{\prime})=n$, then there exist $s_{1}, \cdots, s_{n}$ and $s^{\prime}_{1},\cdots, s^{\prime}_{n}$ such that
	$$g=s_{1}\cdots s_{n}, g^{\prime}=s^{\prime}\cdots s^{\prime}_{n}.
	$$
	
	Take $k=min\{i:s_{i} \neq s^{\prime}_{i}\}$. When $s_{k}<s^{\prime}_{k}$, we have $g
	<g^{\prime}$; otherwise, $g >g^{\prime}$. Then we can arrange the elements in the group. Let $G=\{g_{n}\}^{\infty}_{n=0}$ be an enumeration of $G$ according to the order such that $\ell_{S}(1_{G})=\ell_{S}(g_{0})<\ell_{S}(g_{1})< \ell_{S}(g_{2}) < \ell_{S}(g_{3})\cdots$.
	
	For a finite alphabet $A$, we can define a metric $d$ on $A^{G}$ by 
	$$d(x,y)=2^{-min\left\{\ell_{S}(g_{n})|x_{g_{n}}\neq y_{g_{n}}\right\}}.
	$$
	
	\hspace{4mm}
	For convenience, we always assume $G$ is a finitely generated amenable group in this section and $\{B_{S}(m)\}_{m=1}^{\infty}$ is a F$\phi$lner sequence of $G$.
	
	\hspace{4mm}
	Suppose $A, B$ are two finite sets and $\pi:(A\times B)^{G}\to B^{G}$ is the natural projection. If $((A\times B)^{G}, \sigma)$ and $(B^{G},\sigma)$ are two symbolic dynamical systems, $\Omega $ is a subshift of $(A\times B)^{G}$, and $\Omega^{\prime}=\pi (\Omega)\subseteq B^G$. For any subset $D\subset G$, we define  $\left.\Omega\right|_{D} \subset(A \times B)^{D}$ and  $\left.\Omega^{\prime}\right|_{D} \subset B^{D}$  as the images of the projections of  $\Omega$  and  $\Omega^{\prime}$  to the   $D\subset G$  coordinates, respectively. We denote by  $\pi_{D}:\left.\left.\Omega\right|_{D} \rightarrow \Omega^{\prime}\right|_{D}$  the natural projection map. And we define new metrics  $d $ and  $d^{\prime}$  on  $(A \times B)^{G}$  and  $B^{G}$  by		
	$$d\left((x, y),\left(x^{\prime}, y^{\prime}\right)\right)=2^{-\min \left\{ \ell_{S}(g_{n}) \mid \left(x_{g_{n}}, y_{g_{n}}\right) \neq\left(x_{g_{n}}^{\prime}, y_{g_{n}}^{\prime}\right)\right\}}, \quad d^{\prime}\left(y, y^{\prime}\right)=2^{-\min \left\{\ell_{S}(g_{n}) \mid y_{g_{n}} \neq y_{g_{n}}^{\prime}\right\}}.$$

	\begin{pro}\label{l5.1}
		Let $0 \leq w \leq 1$. If
		$ \pi:\Omega \rightarrow \Omega^{\prime}$ is a factor map, then
		
		$$h_{\text {top }}^{w}(\pi, G )=\lim _{m \rightarrow \infty} \frac{1}{|B_{S}(m)|} \log \sum_{v \in \Omega^{\prime} \mid_{B_{S}(m)}}|\pi_{B_{S}(m)}^{-1}(v)|^{w}.$$ 
		Here $|\pi_{B_{S}(m)}^{-1}(v)|$ is the cardinality of  $\left.\pi_{B_{S}(m)}^{-1}(v) \subset \Omega\right|_{B_{S}(m)}$. 
	\end{pro}
	\begin{proof} 
		Fix $0<\varepsilon<1$ and $m\in \mathbb{N}_0$. For any subsets $U \subset(A \times B)^{G}$ and $V \subset B^{G}$, we denote by  $\left.U\right|_{B_{S}(m)} \subset(A \times B)^{B_{S}(m)}$  and  $\left.V\right|_{B_{S}(m)} \subset B^{B_{S}(m)}$ the projections to the  $B_{S}(m)$ -coordinates. It is easy to see that if  $\operatorname{diam}\left(U, d_{B_{S}(m)}\right)<1$  then  $\left.U\right|_{B_{S}(m)}$  is a singleton or empty. Similarly, if $\operatorname{diam}\left(V, d_{B_{S}(m)}^{\prime}\right)<1$  then so is  $\left.V\right|_{B_{S}(m)}$ .
		Hence	
		$$h_{\text {top }}^{w}(\pi, G) \leq \lim _{m \rightarrow \infty} \frac{1}{|B_{S}(m)|} \log \sum_{\left.v \in \Omega^{\prime}\right|_{B_{S}(m)}}|\pi_{{B_{S}(m)}}^{-1}(v)|^{w}. $$	
		The remaining task is to prove the reverse inequality. 
		
		\hspace{4mm}
		If $\Omega^{\prime}=V_{1}\cup V_{2}\cup \cdots \cup V_{n}$ is an open cover with  $\operatorname{diam}\left(V_{k}, d_{B_{S}(m)}\right)<\varepsilon$ such that	
		$$\#^{w}\left(\pi, \Omega, d_{B_{S}(m)}, d_{B_{S}(m)}^{\prime}, \varepsilon\right)=\sum_{k=1}^{n}\left(\#\left(\pi^{-1}\left(V_{k}\right), d_{B_{S}(m)}, \varepsilon\right)\right)^{w}. $$		
		Assuming that $V_{j}\neq \emptyset$ for any $1\leq j\leq n$. Let $t_{k}:=   \#\left(\pi^{-1}\left(V_{k}\right), d_{B_{S}(m)}, \varepsilon\right)$,
		then we have $	\# \left(\pi, \Omega, d_{B_{S}(m)}, d_{B_{S}(m)}^{\prime}, \varepsilon\right) =\sum_{k=1}^{n} t_{k}^{w}.$
		
		\hspace{4mm}
		For every $1\leq k\leq n$,  we consider an open cover $\{U_{k_ 1}, U_{k_ 2}, \cdots , U_{k_{t_{k}}}\}$ of $\pi^{-1}\left(V_{k}\right)$  with $\text{diam}  \left(U_{k_ l}, d_{B_{S}(m)}\right)<\varepsilon$  for all  $1 \leq l \leq t_{k}$, we also assume that every $U_{k_ l}\neq \emptyset$. So 
		$$
		\sum_{k=1}^{n} t_{k}^{w} =\sum_{\left.v \in \Omega^{\prime}\right|_{B_{S}(m)}}\left(\sum_{k:\left.V_{k}\right|_{B_{S}(m)}=\{v\}} t_{k}^{w}\right).
		$$
		Here, in the expression $\sum_{k: V_{k} \mid _{B_{S}(m)}=\{v\}}$, $ k$  takes over every $k \in[1, n]$ that satisfies  $\left.V_{k}\right|_{B_{S}(m)}=\{v\}.$ Since every $ \left.U_{k _l}\right|_{B_{S}(m)}$  is a singleton, we have
		$$\#(\pi_{B_{S}(m)}^{-1}(v)) \leq \sum_{k:\left.V_{k}\right|_{B_{S}(m)}=\{v\}} t_{k}.$$
		
		\hspace{4mm}
		For each $ \left.v \in \Omega^{\prime}\right|_{B_{S}(m)}$, we have
		$$\pi_{B_{S}(m)}^{-1}(v)=\bigcup_{k:V_{k}\mid _{B_{S}(m)}=\{v\}}(\pi^{-1}(V_{k}))\mid_{B_{S}(m)}=\bigcup_{k: V_{k} \mid _{B_{S}(m)}=\{v\}} \bigcup_{l=1}^{t_{k}}(U_{k l}\mid_{B_{S}(m)}).$$ 
		Recall that $0 \leq w \leq 1$. From the fact that $(p+q)^{w}\leq p^{w}+q^{w}$ for any $p, q\geq 0$, we have 
		$$\sum_{k: V_{k} \mid _{B_{S}(m)}=\{v\}} t_{k}^{w} \geq\left(\sum_{k: V_{k} \mid _{B_{S}(m)}=\{v\}} t_{k}\right)^{w}.$$
		Hence
		\begin{align*}
			\#^{w}\left(\pi, \Omega, d_{B_{S}(m)}, d_{B_{S}(m)}^{\prime}, \varepsilon\right)=& \sum_{\left.v \in \Omega^{\prime}\right|_{B_{S}(m)}}\left(\sum_{k:\left.V_{k}\right|_{B_{S}(m)}=\{v\}} t_{k}^{w}\right) \\
			& \geq \sum_{\left.v \in \Omega^{\prime}\right|_{B_{S}(m)}}\left(\sum_{k:V_{k}|_{B_{S}(m)}=\{v\}} t_{k}\right)^{w} \\
			& \geq \sum_{v \in \Omega^{\prime}|_{B_{S}(m)}}\#(\pi_{B_{S}(m)}^{-1}(v))^{w} . 
		\end{align*}

		Then
		$$ \frac{1}{|B_S(m)|} \log \#^{w}\left(\pi, \Omega, d_{B_{S}(m)}, d_{B_{S}(m)}^{\prime}, \varepsilon\right)\geq  \frac{1}{B_S(m)} \log \sum_{v \in \Omega^{\prime}|_{B_{S}(m)}}\#(\pi^{-1}_{|B_{S}(m)|} (v))^{w}.
		$$
		Letting $m\to \infty$, we conclude 
		\begin{equation*}h_{\text {top }}^{w}(\pi, G) \geq \lim _{m \rightarrow \infty} \frac{1}{|B_S(m)|} \log \sum_{v \in \Omega^{\prime}|_{B_{S}(m)}}\#(\pi^{-1}_{|B_{S}(m)|} (v))^{w}.\qedhere
		\end{equation*}
	\end{proof}

	\subsection{Carpet system}
	\
	\newline \indent \quad\ Now we introduce the infinite-dimensional Bedford-McMullen carpets for amenable group actions. Suppose that $G$ is a finitely generated group with a symmetric generating set $S$ and $a \geq b \geq 2$ are two natural numbers. Write
	$$A=\{0,1,2, \ldots, a-1\}, \quad B=\{0,1,2, \ldots, b-1\} .$$
	Let $(A \times B)^{G}$  be the full $G$-shift over $A \times B$  with the shift action  $\sigma:G\times ((A \times B)^{G}) \rightarrow(A \times B)^{G} $ and $\pi:(A \times B)^{G} \rightarrow B^{G}$  be the natural projection. We also denote by $ \sigma:G\times B^{G} \rightarrow B^{G}$  the full $G$-shift on $ B^{G}$ .
	
	\hspace{4mm}
	Suppose $[0,1]^{G}=\{(x_g)_{g\in G}: x_g\in [0,1]\}$ is an infinite dimensional cube, and $(\alpha_{g})_{g \in G}$ is a sequence of positive real numbers satisfying
	$$
	\alpha_{1_{G}}=1,\sum_{g \in G}\alpha_{g} < +\infty.
	$$
	We define a metric $d$ on the product  $[0,1]^{G}\times[0,1]^{G}$  by
	$$d\left((x, y),\left(x^{\prime}, y^{\prime}\right)\right)=\sum_{g\in G} \alpha_{g} \max \left(\left|x_{g}-x_{g}^{\prime}\right|,\left|y_{g}-y_{g}^{\prime}\right|\right),$$

	where  $(x,y)=(\left(x_{g}\right)_{g \in G}, \left(y_{g}\right)_{g \in G})$ and  $(x^{\prime}, y^{\prime})=(\left(x_{g}^{\prime}\right)_{g \in G}, \left(y_{g}^{\prime}\right)_{g \in G})$  are two points of  $[0,1]^{G}\times[0,1]^{G}.$ We define the shift map $\sigma:G\times([0,1]^{G} \times[0,1]^{G}) \rightarrow[0,1]^{G} \times[0,1]^{G}$   by
	$$h\left(\left(x_{g}\right)_{g \in G},\left(y_{g}\right)_{g \in G}\right)=\left(\left(x_{gh}\right)_{g \in G},\left(y_{gh}\right)_{g \in G}\right), \text{ for all } h\in G .$$
	The shift $\sigma_1$ and $\sigma_2$ on $[0,1]^{G} \times[0,1]^{G}$ is defined by $((\sigma_{1,g}x)_{(g_1,g_2)}, (\sigma_{1,g}y)_{(g_1,g_2)})=(x_{(g_1g,g_2)}, y_{(g_1g,g_2)})$ and $((\sigma_{2,h}x)_{(g_1,g_2)}, (\sigma_{2,h}y)_{(g_1,g_2)})=(x_{(g_1,g_2h)}, y_{(g_1,g_2h)})$ for every $g\in G_1, h\in G_2, (g_1, g_2)\in G$. 
	
	\hspace{4mm}
	Let  $\Omega \subset(A \times B)^{G} $ be a subshift. We assume  $\Omega \neq \emptyset $. We define a carpet system $X_{\Omega} \subset[0,1]^{G} \times[0,1]^{G}$  by
	$$X_{\Omega}=\left\{\left(\sum_{m=1}^\infty \frac{x_{m}}{a^{m}}, \sum_{m=1}^\infty \frac{y_{m}}{b^{m}}\right) \in[0,1]^{G} \times[0,1]^{G} \mid\left(x_{m}, y_{m}\right) \in \Omega \text { for any } m\geq 1\right\} .$$
	Here $ x_{m} \in A^{G} \subset \ell^{\infty}$  and   $\sum_{m=1}^\infty \frac{x_{m}}{a^{m}}  \in  \ell^{\infty}$, then we have $\sum_{m=1}^\infty \frac{x_{m}}{a^{m}} \in   [0,1]^{G}$ . Similarly for the term  $\sum_{m=1}^\infty \frac{y_{m}}{b^{m}} $. Therefore $\left(X_{\Omega}, \sigma\right)$  is a subsystem of  $\left([0,1]^{G} \times[0,1]^{G}, \sigma\right).$ Define $ \Omega^{\prime}=\pi(\Omega) $, then $\Omega^{\prime}$ is a subshift of  $B^{G}$.

	\subsection{The calculation for carpet systems}
	\
	\newline 
	\indent\quad\quad	In the following, we are going to compute the metric mean dimension and mean Hausdorff dimension of the infinite-dimensional Bedford-McMullen carpets, respectively. Before our calculation, we need some preparations. Set
	$w=\log _{a} b.$ As $a \geq b \geq 2$ , we have  $0<w \leq 1$. For any natural number  $m$, we set
	
	$$\left.\Omega^{\prime}\right|_{B_{S}(m)}=\{\left(v_{n}\right)_{n \in B_S(m)}\in B^{B_S(m)}\mid \left(v_{n}\right)_{n \in G}\in B^{G}\},$$
	\begin{align*}
		&\left.\Omega\right|_{B_{S}(m)}\\
		&=\{\left(\left(u_{n}\right)_{n \in B_S(m)},\left(v_{n}\right)_{n \in B_S(m)}\right)\in A^{B_S(m)} \times B^{B_S(m)}\mid \left(\left(u_{n}\right)_{n \in G},\left(v_{n}\right)_{n \in G}\right)\in (A \times B)^G\},\\
		&\left.X_{\Omega}\right|_{ B_S(m)}\\
		&=\left\{\left(\sum_{n=1}^\infty \frac{x_{n}}{a^{n}}, \sum_{n=1}^\infty \frac{y_{n}}{b^{n}}\right) \in[0,1]^{ B_S(m)} \times[0,1]^{ B_S(m)}: \left(x_{n}, y_{n}\right) \in \Omega|_{ B_S(m)} \text { for all } n\geq 1\right\} .
	\end{align*}
	
	\hspace{4mm}
	Let  $(x, y) \in\left(\left.\Omega\right|_{B_S(m)}\right)^{\mathbb{N}} $ where  $x=\left(x_{m}\right)_{m=1}^\infty$  and  $y=\left(y_{m}\right)_{m=1}^\infty $ with  $x_{m} \in A^{ B_S(m)}, y_{m} \in B^{ B_S(m)}$  and  $\left.\left(x_{m}, y_{m}\right) \in \Omega\right|_{ B_S(m)} .$ For natural numbers  $m$  and  $l$ , we define a subset  $\left.\Phi_{l, m}(x, y) \subset X_{\Omega}\right|_{ B_S(m)}$  by
	
	$$\Phi_{l, m}(x, y)=\left\{\left(\sum_{n=1}^\infty \frac{x_{n}^{\prime}}{a^{n}}, \sum_{n=1}^\infty \frac{y_{n}^{\prime}}{b^{n}}\right) \mid \begin{array}{c}
		\left.\left(x_{n}^{\prime}, y_{n}^{\prime}\right) \in \Omega\right|_{ B_S(m)} \text { for all } n\geq 1 \text { with } \\
		x_{n}^{\prime}=x_{n}(1\leq n\leq \lfloor wl\rfloor) \text { and } \\
		y_{n}^{\prime}=y_{n}(1\leq n\leq l)
	\end{array}\right\} .$$

	Here  $\lfloor wl\rfloor $ is the largest integer not greater than  $wl $. The set $ \Phi_{l, m}(x, y)$  depends only on the coordinates  $(x_{i})_{i=1}^{\lfloor w l\rfloor}, (y_{i})_{i=1}^l.$ So we also denote it by
	
	$$\Phi_{l, m}\left((x_{i})_{i=1}^{\lfloor w l\rfloor}, (y_{i})_{i=1}^l\right) \quad\left(=\Phi_{l, m}(x, y)\right) .$$
	We also define a subset $\Psi_{l, m}(x, y)(\Omega|_{ B_S(m)})^{\mathbb{N}}$ by
	
	$$\Psi_{l, m}(x, y)=\left\{\left(x^{\prime}, y^{\prime}\right) \in\left(\left.\Omega\right|_{ B_S(m)}\right)^{\mathbb{N}} \mid \begin{array}{c}
		x_{n}^{\prime}=x_{n}( 1\leq n\leq \lfloor wl\rfloor) \\
		y_{n}^{\prime}=y_{n}(1\leq n\leq l)
	\end{array}\right\} .$$
	then the set  $\Phi_{l, n}(x, y) $  is the image of  $\Psi _{l, n}(x, y)$  under the map $$(\Omega\mid _{B_S(m)})^{\mathbb{N}} \rightarrow X_{\Omega}|_{B_S(m)}, \quad (x^{\prime}, y^{\prime}) \mapsto (\sum_{n=1}^\infty \frac{x_{n}^{\prime}}{a^{n}}, \sum_{n=1}^\infty \frac{y_{n}^{\prime}}{b^{n}}) .$$
	It is straightforward to check that for  $(x, y)=\left(x_{n}, y_{n}\right)_{n \in \mathbb{N}}$  and  $\left(x^{\prime}, y^{\prime}\right)=\left(x_{n}^{\prime}, y_{n}^{\prime}\right)_{n \in \mathbb{N}}$  in  $\left(\left.\Omega\right|_{B_S(m)}\right)^{\mathbb{N}} $, we have the following relationship:\\

	\begin{tikzcd}[row sep=tiny]
		& \Psi_{l, m}(x, y)= \Psi_{l, m}\left(x^{\prime}, y^{\prime}\right) \arrow[dd,Leftrightarrow] \\
		\left(x_{1}, \ldots, x_{\lfloor w l\rfloor}, y_{1}, \ldots, y_{l}\right)=\left(x_{1}^{\prime}, \ldots, x_{\lfloor w l\rfloor}^{\prime}, y_{1}^{\prime}, \ldots, y_{l}^{\prime}\right) \arrow[ur,Leftrightarrow] \arrow[dr,Leftrightarrow] & \\
		& \Phi _{l, m}(x, y)= \Phi _{l, m}\left(x^{\prime}, y^{\prime}\right)
	\end{tikzcd}

	\hspace{4mm}
	For any subsets  $E$  and  $F$  of  $\left.X_{\Omega}\right|_{B_S(m)},$ we set
	
	$$\operatorname{dist}_{\infty}(E, F)=\inf _{\substack{x \in E \\ y \in F}}\|x-y\|_{\infty} .$$

	Here  $\|x-y\|_{\infty}$  is the  $\ell^{\infty}$ -distance on  $\left.X_{\Omega}\right|_{B_S(m)} \subset (\mathbb{R}^2)^{B_S(m)}$ .
	
	\hspace{4mm}
	The following lemma shows the distance of each set $\Phi_{l,m}(x,y)$ under a special condition.
	
	\begin{lem}\label{l5.13}
		Let  $l$  and  $m$  be natural numbers. Let  $\left(x^{(1)}, y^{(1)}\right), \ldots,\left(x^{(k)}, y^{(k)}\right) \in\left(\left.\Omega\right|_{B_S(m)}\right)^{\N}$  with  $k \geq 4^{|B_S(m)| }+1$ . Suppose that
		$$\Psi_{l, m}\left(x^{(i)}, y^{(i)}\right) \neq \Psi_{l, m}\left(x^{(j)}, y^{(j)}\right) \quad \text { for } i \neq j .$$
		Then there are  $i$  and  $j$  for which
		$$\operatorname{dist}_{\infty}\left(\Phi_{l, m}\left(x^{(i)}, y^{(i)}\right), \Phi_{l, m}\left(x^{(j)}, y^{(j)}\right)\right) \geq b^{-L} .$$
	\end{lem}
	\begin{proof}
		The assumption  $\Psi_{l, m}\left(x^{(i)}, y^{(i)}\right) \neq \Psi_{l, m}\left(x^{(j)}, y^{(j)}\right) $ implies
		
		$$\left(x_{1}^{(i)}, \ldots, x_{\lfloor w l\rfloor}^{(i)}, y_{1}^{(i)}, \ldots, y_{l}^{(i)}\right) \neq\left(x_{1}^{(j)}, \ldots, x_{\lfloor w l\rfloor}^{(j)}, y_{1}^{(j)}, \ldots, y_{l}^{(j)}\right) .$$

		Since  $k \geq 4^{|B_S(m)|}+1 $, there are  $i$  and  $j$  for which we have either
		
		$$\left\|\sum_{n=1}^{\lfloor w l\rfloor} \frac{x_{n}^{(i)}}{a^{n}}-\sum_{n=1}^{\lfloor w l\rfloor} \frac{x_{n}^{(j)}}{a^{n}}\right\|_{\infty} \geq 2 a^{-\lfloor w l\rfloor}$$
		
		or
		
		$$\left\|\sum_{n=1}^{l} \frac{y_{n}^{(i)}}{b^{n}}-\sum_{n=1}^{l} \frac{y_{n}^{(j)}}{b^{n}}\right\|_{\infty} \geq 2 b^{-l}.$$

		We have
		
		$$\Phi_{l, m}(x, y) \subset\left(\sum_{n=1}^{\lfloor w l\rfloor} \frac{x_{n}}{a^{n}}, \sum_{n=1}^{m} \frac{y_{n}}{b^{n}}\right)+\left[0, a^{-\lfloor w l\rfloor}\right]^{|B_S(m)|} \times\left[0, b^{-l}\right]^{|B_S(m)| } .$$

		and either
		
		$$\operatorname{dist}_{\infty}\left(\Phi_{l, m}\left(x^{(i)}, y^{(i)}\right), \Phi_{l, m}\left(x^{(j)}, y^{(j)}\right)\right) \geq a^{-\lfloor w l\rfloor} \geq b^{-l} ,$$

		or
		
		$$\operatorname{dist}_{\infty}\left(\Phi_{l, m}\left(x^{(i)}, y^{(i)}\right), \Phi_{l, m}\left(x^{(j)}, y^{(j)}\right)\right) \geq b^{-l} .$$

		So we get the statement in both of the cases.
	\end{proof}

	\subsubsection{Mean Hausdorff dimension of the carpet system}
	\
	\newline \indent \quad\ For  $\left.v \in \Omega^{\prime}\right|_{B_S(m)}$, we set  
	\begin{align*}
		t_{m}(v)&=\#\{u \in A^{ B_S(m)}\mid \left.(u, v) \in \Omega\right|_{ B_S(m)}\},\\
		Z_{m}&=\sum_{\left.v \in \Omega^{\prime}\right|_{ B_S(m)}} t_{m}(v)^{w}.
	\end{align*}
	Then
	$$h_{\text {top }}^{w}(\pi, G)=\lim _{m \rightarrow \infty} \frac{\log Z_{m}}{B_S(m)} \qquad (\text{ by Lemma } \ref{l5.1}).$$
	\hspace{4mm}
	For any $\left.(u, v) \in \Omega\right|_{B_S(m)} $, we define $f_{m}(u, v)$ as 
	$$f_{m}(u, v)=\frac{1}{Z_{m}} t_{m}(v)^{w-1} .$$
	Notice that 
	$$\sum_{\left.(u, v) \in \Omega\right|_{ B_S(m)}} f_{m}(u, v)=1 .$$
	Thus $f_{m}(u, v)$ can be regarded as a probability measure $\mu_{m}$ on  $\left.\Omega\right|_{B_S(m)}$. We can also define a probability measure on  $\left(\left.\Omega\right|_{ B_S(m)}\right)^{\N}$ by
	$\mu_{m}=\left(f_{m}\right)^{\otimes \mathbb{N}}$. It is the product of infinite copies of the measure $f_{m}$. 
	
	\hspace{4mm} 
	The following lemmas are crucial, the first one comes directly from \cite{TM}, while the second one is similar to Lemma 5.11 in \cite{TM} and its proof follows the same approach, so we omit it.
	
	\begin{lem}\label{l57}
		Let  $c, \varepsilon, s $ be positive numbers. Let  $(X, d) $ be a compact metric space with a Borel probability measure  $\mu$ . If $6 \varepsilon^{c}<1 $ and for any  $x \in X,$ there exists a Borel subset  $A \subset X $ containing  $x$  and
		$$0<\operatorname{diam} A<\frac{\varepsilon}{6}, \quad \mu(A) \geq(\operatorname{diam} A)^{s},$$
		then  $\operatorname{dim}_{\mathrm{H}}(X, d, \varepsilon) \leq(1+c) s .$
	\end{lem}
	
	\begin{lem} \label{l5.11}
		For any $\delta>0$ and any $K\in \N$,  there exists a natural number  $N=N(\delta, K) \geq K$  for which the following statement holds true: For any  $m \in \N_0$  and any  $(x, y) \in\left(\left.\Omega\right|_{ B_S(m)}\right)^{\N} $, there exists a natural number  $l \in[K, N]$  satisfying
		$$\frac{1}{|B_S(m)| l} \log \mu_{m}\left(\Psi_{l, m}(x, y)\right) \geq-\frac{\log Z_{m}}{| B_S(m)|}-\delta .$$
	\end{lem}
	
	\hspace{4mm}
	Next, we will first establish an upper bound for the Hausdorff dimension of $X_{\Omega}$,  based on the idea presented in Lemma \ref{l57}. From this, we will derive the upper bound for the mean Hausdorff dimension.

	\begin{thm}\label{t5.16}
		For the carpet system $X_{\Omega}$, we have $$\overline{\operatorname{mdim}}_{\mathrm{H}}\left(X_{\Omega}, \{B_S(m)\}, d\right) \leq \frac{h_{\text {top }}^{w}(\pi, G)}{\log b} .$$
	\end{thm}
	\begin{proof} As the map $B_S(m)\to Z_m$ satisfies the conditions in Lemma \ref{l100}, we can find a $t>0$ such that $\log _{b} Z_{m} \leq t \cdot | B_S(m)|$ holds for all  $m \geq 1$.

		\hspace{4mm}
		Fix  $m\in \N$. We define $\nu_{m}$ as the push-forward of  $\mu_{m}$  under the map
		
		$$(\Omega|_{B_S(m)})^{\mathbb{N}} \rightarrow X_{\Omega}|_{B_S(m)}, \quad(x_{n}, y_{n})_{n \in \mathbb{N}} \mapsto(\sum_{n=1}^{\infty} \frac{x_{n}}{a^{n}}, \sum_{n=1}^{\infty} \frac{y_{n}}{b^{n}}) .$$
		It is a probability measure on  $\left.X_{\Omega}\right|_{ B_S(m)}$. And we define $X(m)$ by
		$$X(m)=X_{\Omega}|_{ B_S(m)} \times[0,1]^{2} \times[0,1]^{2} \times \cdots \subset([0,1]^{2})^{G}.$$
		As $X_{\Omega}|_{B_S(m)} \subset\left([0,1]^{2}\right)^{ B_S(m)})$, we have $X_{\Omega} \subset X(m) $. 
		
		\hspace{4mm}
		Denote by $Leb$ the Lebesgue measure on the square  $[0,1]^{2}$. Then we can define a probability measure $\nu_{m}^{\prime}$ on  $X(m)$ by
		$$\nu_{m}^{\prime}:=\nu_{m} \otimes L e b \otimes L e b \otimes L e b \otimes \cdots .$$
		For any $L, m\in \N $ and any $(x, y) \in\left(\left.\Omega\right|_{ B_S(m)}\right)^{\mathbb{N}}$, write
		$$\Phi_{l, m}^{\prime}(x, y):= \Phi_{l, m}(x, y) \times[0,1]^{2} \times[0,1]^{2} \times \cdots.$$
		Then we have $\Phi_{l, m}^{\prime}(x, y)\subset X(m) $ and 
		$$\nu_{m}^{\prime}\left(\Phi_{l, m}^{\prime}(x, y)\right)=\nu_{m}\left(\Phi_{l, m}(x, y)\right) \geq \mu_{m}\left(\Psi_{l, m}(x, y)\right) .$$
		
		\hspace{4mm}
		Fix  $\varepsilon>0 $, we choose  $\delta>0$  with  $6 \delta^{\varepsilon}<1 .$ Let $R=\sum_{g\in \Z^k}\alpha_g$, we take $K\in \N$  satisfying	
		$$a b^{-K}<\frac{\delta}{6R}, \quad(Ra)^{s+\varepsilon} \leq b^{\frac{\varepsilon K}{2}} .$$	
		Let  $N=N\left(\frac{\varepsilon \log b}{2}, K\right) $ be the natural number given by Lemma \ref{l5.11}, then by Lemma \ref{l5.11} for any  $m \geq 1$  and any  $(x, y) \in\left(\left.\Omega\right|_{B_S(m)}\right)^{\mathbb{N}} $, there exists a natural number $l \in[K, N] $ such that
		
		$$\frac{1}{l \cdot |B_S(m)|} \log \mu_{m}\left(\Psi_{l, m}(x, y)\right) \geq-\frac{\log Z_{m}}{| B_S(m)|}-\frac{\varepsilon \log b}{2} .$$
		
		That is
		$$\mu_{m}(\Psi_{l, m}(x, y))\geq \exp \left\{-l \log b\left(\log _{b} Z_{m}+\frac{\varepsilon \cdot | B_S(m)|}{2}\right)\right\}.$$
		This gives that $$\nu_{m}^{\prime}(\Phi_{l, m}^{\prime}(x, y))\geq \left(b^{-l}\right)^{\log _{b} Z_{m}+\frac{\varepsilon \cdot |B_S(m)|}{2}}.$$

		\begin{claim} \label{c5.17}
			For each  $m \geq 1$ and  $(x, y) \in\left(\Omega|_{ B_S(m)}\right)^{\mathbb{N}}$, there exists a natural number $l \in[K, N]$  such that
			$$\nu_{m}^{\prime}\left(\Phi_{l, m}^{\prime}(x, y)\right) \geq\left(Ra b^{-l}\right)^{\log _{b} Z_{m}+| B_S(m)|\varepsilon} .$$
		\end{claim}
		
		\begin{proof} By the above argument, we can find a positive integer  $l \in[K, N] $ such that
			$$\nu_{m}^{\prime}\left(\Phi_{l, m}^{\prime}(x, y)\right) \geq\left(b^{-l}\right)^{\log _{b} Z_{m}+\frac{| B_S(m)|\varepsilon }{2}} .$$
			As $\log Z_{m} \leq | B_S(m)|\cdot s $ for any $m\in \N_0$ and $l \geq K$, then
			$$(R a)^{\frac{\log _{b} Z_{m}}{| B_S(m)|}+\varepsilon} \leq(R a)^{s+\varepsilon} \leq b^{\frac{\varepsilon K}{2}}\leq b^{\frac{\varepsilon l}{2}}.$$
			Thus
			$$\left(b^{-l}\right)^{\log _{b} Z_{m}+\frac{| B_S(m)|\varepsilon }{2}} \geq\left(R a b^{-l}\right)^{\log _{b} Z_{m}+\varepsilon |B_S(m)|} .$$
			This gives the statement.		
		\end{proof}
		
		\hspace{4mm}
		Now we pick a $ m_{0}\in \N$  such that
		$$\sum_{g\in G\setminus B_S(m_{0})} \alpha_g<a b^{-N},$$
		and let $c=\sum_{g\in B_S(m_0) } \alpha_g$.  For any $u, v \in\left([0,1]^{2}\right)^{G}$ with $u=(p_g, q_g)_{g\in G}, v=(p_g^{\prime}, q_g^\prime)_{g\in G})$  and  any natural number $m$ with $m>m_0$, we have	
		\begin{align*}
			d_{B_S(m-m_0)}(u, v)&=\max_{i\in B_S(m-m_0) } \sum_{g\in G} \alpha_g \max \left(\left|p_{gi}-p_{gi}^{\prime}\right|,\left|q_{gi}-q_{gi}^{\prime}\right|\right)\\
			&=\max_{i\in B_S(m-m_0)}\{\sum_{g\in B_S(m_0)} \alpha_g \max \left(\left|p_{gi}-p_{gi}^{\prime}\right|,\left|q_{gi}-q_{gi}^{\prime}\right|\right)\\
			&\quad +\sum_{g\in G\setminus B_S(m_0)} \alpha_g \max \left(\left|p_{gi}-p_{gi}^{\prime}\right|,\left|q_{gi}-q_{gi}^{\prime}\right|\right)\}\\
			&\leq \max_{g\in B_S(m)}\left(|p_{g}-p_{g}^{\prime}|,|q_{g}-q_{g}^{\prime}|\right)\cdot \sum_{g\in B_S(m_0)} \alpha_g+\sum_{g\in G\setminus B_S(m_0)} \alpha_g\\
			&\leq c\left\|\left.u\right|_{B_S(m)}-\left.v\right|_{ B_S(m)}\right\|_{\infty}+a b^{-N} .
		\end{align*}
		
		Recall that $$\Phi_{l, m}^{\prime}(x, y)= \Phi_{l, m}(x, y) \times[0,1]^{2} \times[0,1]^{2} \times \cdots .$$
		Then 
		$$0 \leq \operatorname{diam}\left(\Phi_{l, m}^{\prime}(x, y), d_{ B_S(m-m_0)}\right)\leq c\operatorname{diam}\left(\Phi_{l, m}(x, y),\|\cdot\|_{\infty}\right)+a b^{-N} .$$
		
		\begin{claim} \label{c5.18}
			For any $(x, y) \in\left(\left.\Omega\right|_{ B_S(m)}\right)^{\mathbb{N}}$, any $l, m\in \N$ with  $m>m_{0}$  and  $l \leq N $, 
			$$0<\operatorname{diam}\left(\Phi_{L, m}^{\prime}(x, y), d_{ B_S(m-m_0)}\right)<\delta / 6.$$	
		\end{claim}
		
		\begin{proof} 
			Recall that  $\operatorname{diam}\left(\Phi_{l, m}(x, y),\|\cdot\|_{\infty}\right) \leq a b^{-l} $, combine with $c+1< R$ and $l \leq N$, we get
			\begin{equation*}\operatorname{diam}\left(\Phi_{l, m}^{\prime}(x, y), d_{ B_S(m-m_0)}\right)\leq c a b^{-l}+a b^{-N} < R a b^{-l}<\delta / 6.\qedhere
			\end{equation*}
		\end{proof}

		\hspace{4mm}
		Recall that $6 \delta^{\varepsilon}<1$ and $\Phi_{l, m}^{\prime}(x, y)\in X(m)$. Combining with Claims \ref{c5.17} and \ref{c5.18}, for any natural number $m>m_{0}$,
		$$\operatorname{dim}_{\mathrm{H}}\left(X(m), d_{ B_S(m-m_0)}, \delta\right) \leq(1+\varepsilon)\left(\log _{b} Z_{m}+|B_S(m)|\varepsilon \right), \quad ( \text{ By Lemma } \ref{l57}.)$$
		As $X_{\Omega} \subset X(m)$ and $0<\delta<\varepsilon$, we have   $$\operatorname{dim}_{\mathrm{H}}\left(X_{\Omega}, d_{ B_S(m-m_0)}, \varepsilon\right) \leq(1+\varepsilon)\left(\log _{b} Z_{m}+| B_S(m)|\varepsilon \right), \quad\left(m>m_{0}\right) .$$
		Then
		$$\frac{\operatorname{dim}_{\mathrm{H}}\left(X_{\Omega}, d_{ B_S(m-m_0)}, \varepsilon\right)}{|B_S(m)|} \leq(1+\varepsilon)\left(\frac{\log _{b} Z_{m}}{| B_S(m)|}+\varepsilon\right) .$$	
		Recall that 
		$$h_{\text {top }}^{w}(\pi, G)=\lim _{m \rightarrow \infty} \frac{\log  Z_{m}}{| B_S(m)|}.$$ Therefore 
		\begin{align*}
			\limsup _{m \rightarrow \infty} \frac{\operatorname{dim}_{\mathrm{H}}\left(X_{\Omega}, d_{B_S(m)}, \varepsilon\right)}{| B_S(m)|} & \leq (1+\varepsilon)\left(\lim _{m \rightarrow \infty} \frac{\log _{b} Z_{m}}{| B_S(m)|}+\varepsilon\right) \\
			& =\lim_{\varepsilon \rightarrow 0}\frac{1+\varepsilon}{\log b}\left(h_{\mathrm{top}}^{w}(\pi, G)+\varepsilon\right) .
		\end{align*}
		Letting $\varepsilon\to 0$, we obtain the desired result.
	\end{proof}

	\hspace{4mm}
	Next, we prove the lower bound for the mean Hausdorff dimension. We need the following lemma, which is as analogous to that of Tsukamoto \cite{TM}, we omit the proof here.
	
	\begin{lem}\label{l5.12}
		For any $ \delta>0$ there exists a natural number  $m_{2}=m_{2}(\delta)>0$  such that for any natural number $m$, there is a Borel subset  $R(\delta, m) \subset\left(\left.\Omega\right|_{ B_S(m)}\right)^{\mathbb{N}} $ with $ \mu_{m}(R(\delta, m)) \geq \frac{1}{2} $, and for any $(x, y) \in R(\delta, m)$ and $l \geq m_{2}$,  	
		$$\left|\frac{1}{l\cdot |B_S(m)|} \log \mu_{m}\left(\Psi_{l, m}(x, y)\right)+\frac{\log Z_{m}}{|B_S(m)|}\right| \leq \delta .$$
	\end{lem}

	\hspace{4mm}
	With these lemmas, we give a lower bound for the Hausdorff dimension of $X_{\Omega}$.
	\begin{pro} \label{p5.14}
		For any  $\delta>0 $, there is  $\varepsilon>0$  such that for any  $m\in \N$, we have	$$\operatorname{dim}_{\mathrm{H}}\left(\left.X_{\Omega}\right|_{B_S(m)},\|\cdot\|_{\infty}, \varepsilon\right) \geq \log _{b} Z_{m}-\delta \cdot |B_S(m)| .$$
	\end{pro}
	
	\begin{proof}
		We may assume that $\#(\Omega|_{B_S(m)})\geq 2$, for otherwise the assertion is trivial. As in the proof of Proposition \ref{t5.16}, we pick a $s>0$  satisfying $ \log _{b} Z_{m} \leq   s\cdot |B_S(m)|$  for all  $m\in \N_0$ . Fix a $\delta>0,$ let $m_{2}$ be as in Lemma \ref{l5.12} for $(\frac{\delta \log b}{2})$. By Lemma \ref{l5.12}, for any $m\in \N_0,$  there is a Borel subset $ R=R\left(\frac{\delta \log b}{2}, m\right) \subset\left(\left.\Omega\right|_{B_S(m)}\right)^{\N}$  satisfying  $\mu_{m}(R) \geq 1 / 2$  and for any  $l \geq m_{2}$  and  $(x, y) \in\left(\left.\Omega\right|_{B_S(m)}\right)^{\N} $, 
		$$\left|\frac{1}{ l\cdot |B_S(m)|} \log \mu_{m}\left(\Psi_{l, m}(x, y)\right)+\frac{\log Z_{m}}{|B_S(m)|}\right| \leq \frac{\delta \log b}{2} .$$
		Fix a $m\in \N_0$ and pick a $\varepsilon>0$  such that
		$$\varepsilon<b^{-m_{2}-1}, \quad 4 b^{\varepsilon} \varepsilon^{\delta / 2}<\frac{1}{2} .$$
		
		\hspace{4mm}
		If $C$ is a subset of $X$, we define $D(C):=\operatorname{diam}\left(C,\|\cdot\|_{\infty}\right)$ for simplicity. Let $\{E_1, E_2, E_3, \cdots\}$ be a cover of $\left.X_{\Omega}\right|_{B_S(m)}$ with $0<D\left(E_{k}\right)<\varepsilon $ for every $ k \geq 1 $. For every $E_{k}$, we pick a natural number $l_k$ with $l_{k} \geq m_{2}$ such that
		$$b^{-l_{k}-1} \leq D\left(E_{k}\right)<b^{-l_{k}} .$$
		For every $(x, y)$ in $R$ with $\Phi_{ l_{k}, m}(x, y) \cap E_{k} \neq \emptyset$, there is a unique $\Psi_{l_{k}, m}(x, y)$, we denote by $\mathcal{G}_{k}$ the collection of all such $\Psi_{l_{k}, m}(x, y)$. Then $$\mu_{m}(P) \leq\left(b \cdot D\left(E_{k}\right)\right)^{\log _{b} Z_{m}-\frac{\delta  \cdot |B_S(m)|}{2}}, \quad\left(P \in \mathcal{C}_{k}\right) .$$
		
		\hspace{4mm}
		For every $l_{k} \geq m_{2} $ and every  $Q=\Psi_{l_{k}, m}(x, y) \in \mathcal{G}_{k}$  with  $(x, y) \in R$,

		$$\frac{\delta \log b}{2}-\frac{\log Z_{m}}{|B_S(m)|} \geq \frac{1}{l_{k}\cdot |B_S(m)|} \log \mu_{m}(Q).$$
		That is
		\begin{align*}
			b^{-l_{k}\left(\log _{b} Z_{m}-\frac{\delta \cdot |B_S(m)|}{2}\right)}\geq \mu_{m}(P) .
		\end{align*}
		For any $P, Q\in \mathcal{G}_{k}$, we have	$$\operatorname{dist}_{\infty}\left(P,Q\right) \leq D\left(E_{k}\right)<\varepsilon<b^{-l_{k}}.$$
		by Lemma \ref{l5.13}, $\left|\mathcal{G}_{k}\right| \leq 4^{B_S(m)}$. As $$R\subset\bigcup_{k=1}^\infty\bigcup_{P\in\mathcal{G}_k}P,$$
		then $$\frac{1}{2}\leq\mu_m(R)\leq\sum_{k=1}^\infty\sum_{P\in\mathcal{G}_k}\mu_m(P),$$
		and
		\begin{align*}
			\sum_{k=1}^\infty\sum_{P\in\mathcal{G}_k}\mu_m(P)&\leq \sum_{k=1}^{\infty} 4^{|B_S(m)|}\left(b \cdot D\left(E_{k}\right)\right)^{\log _{b} Z_{m}-\frac{\delta\cdot |B_S(m)|}{2}}\\
			&\leq \sum_{k=1}^{\infty} 4^{|B_S(m)|} b^{\log _{b} Z_{M}} D\left(E_{k}\right)^{\log _{b} Z_{m}-\frac{\delta \cdot |B_S(m)|}{2}} \quad(\text{by } b^{-l_{k}-1} \leq D\left(E_{k}\right))\\
			&\leq \sum_{k=1}^{\infty} 4^{|B_S(m)|} b^{s \cdot |B_S(m)|} D(E_{k})^{\log _{b} Z_{m}-\frac{\delta \cdot |B_S(m)|}{2}} \quad(\text{by } \log _{b} Z_{m} \leq   s\cdot |B_S(m)|).
		\end{align*}
		Recall that $D\left(E_{k}\right)<\varepsilon$ and $4 b^{s} \varepsilon^{\delta / 2}<\frac{1}{2},$ we have
		$$\sum_{k=1}^{\infty} \frac{1}{2^{|B_S(m)|}} \cdot D\left(E_{k}\right)^{\log _{b} Z_{m}-\delta\cdot |B_S(m)|}\geq \frac{1}{2}.$$
		Hence
		$$\sum_{k=1}^{\infty} D\left(E_{k}\right)^{\log _{b} Z_{m}-\delta \cdot | B_S(m)|}>1 .$$
		This shows $\operatorname{dim}_{\mathrm{H}}\left(\left.X_{\Omega}\right|_{B_S(m)},\|\cdot\|_{\infty}, \varepsilon\right) \geq \log _{b} Z_{m}-\delta \cdot |B_S(m)| .$
	\end{proof}
	
	\hspace{4mm}
	Now we give the lower bound for the mean Hausdorff dimension for $ X_{\Omega}$.
	
	\begin{cor} \label{c5.15}
		
		$$\underline{\operatorname{mdim}}_{\mathrm{H}}\left(X_{\Omega}, \{B_S(m)\}, d\right) \geq \frac{h_{\text {top }}^{w}(\pi, G)}{\log b} .$$
	\end{cor}
	
	\begin{proof}
		
		For each $m\in \N$ , the natural projection  $\left(X_{\Omega}, d_{B_S(m)}\right) \rightarrow\left(\left.X_{\Omega}\right|_{B_S(m)},\|\cdot\|_{\infty}\right)$  is one-Lipschitz. Thus for each  $\varepsilon>0$, we have
		$$\operatorname{dim}_{\mathrm{H}}\left(X_{\Omega}, d_{B_S(m)}, \varepsilon\right) \geq \operatorname{dim}_{\mathrm{H}}\left(\left.X_{\Omega}\right|_{B_S(m)},  \|\cdot\|_{\infty}, \varepsilon\right) .$$
		By Proposition \ref{p5.14}, for any  $\delta>0$  we can find a  $\varepsilon>0$  such that for any  $m \geq 1$,	$$\operatorname{dim}_{\mathrm{H}}\left(X_{\Omega}, d_{B_S(m)}, \varepsilon\right) \geq \operatorname{dim}_{\mathrm{H}}\left(\left.X_{\Omega}\right|_{B_S(m)},  \|\cdot\|_{\infty}, \varepsilon\right) \geq \log _{b} Z_{m}-\delta \cdot |B_S(m)| .$$
		Then 
		\begin{align*}
			\lim_{\v\to 0}\liminf _{m \rightarrow \infty} \frac{\operatorname{dim}_{\mathrm{H}}\left(X_{\Omega}, d_{B_S(m)}, \varepsilon\right)}{| B_S(m)|} & \geq \lim_{\v\to 0}\lim _{m \rightarrow \infty} \frac{\log _{b} Z_{m}}{|B_S(m)| }-\delta\\
			&=\lim_{\v\to 0}\frac{1}{\log b} \lim _{m \rightarrow \infty} \frac{\log Z_{m}}{| B_S(m)|}-\delta \\
			&=\lim_{\v\to 0}\frac{h_{\text {top }}^{w}(\pi, G)}{\log b}-\delta.
		\end{align*}

		Letting  $\delta \rightarrow 0$ , we obtain the statement.
	\end{proof} 
	
	\hspace{4mm}
	With Corollaries \ref{t5.16} and \ref{c5.15}, we calculate the mean Hausdorff dimension of the carpet system $\left(X_{\Omega}, G\right)$.
	
	\begin{cor}
		The mean Hausdorff dimension of the carpet system $X_{\Omega}$ is given by 
		$$\operatorname{mdim}_{\mathrm{H}}\left(X_{\Omega}, \{B_S(m)\},  d\right) =\frac{h_{\text {top }}^{w}(\pi, G)}{\log b} .$$
	\end{cor}
		
	\subsubsection{Metric mean dimension of the carpet system}	
	\
	\newline \indent \quad\ We will use $\epsilon$-covering number as main tool to begin our calculation. Firstly, we recall the notion of separated set.

 \hspace{4mm}
Let $E$ be a compact subset of $X$, $F \in Fin(G)$ and $\varepsilon>0$. $K \sub  E $ is an $(F,\varepsilon)$-separated set of $E$ if we have
$d_{F}(x,y) \geq \varepsilon$ for any distinct $x,y\in K$.
	\begin{lem}\label{l5.6}
		Let $\|\cdot\|_{\infty}$ be the distance defined by the  $\ell^{\infty}$ -norm on  $\left.X_{\Omega}\right|_{B_S(m)} \subset (\mathbb{R}^2)^{ B_S(m)}$. For any $ l,m\in \N$,  we have	
		\begin{align*}
			\#(X_{\Omega}|_{B_S(m)},\|\cdot\|_{\infty}, b^{-l}) \geq&\#(\Omega|_{B_S(m)})^{\lfloor w l\rfloor} \cdot\#(\Omega^{\prime}|_{B_S(m)})^{l-\lfloor w l\rfloor},\\
			\#(X_{\Omega}|_{B_S(m)},\|\cdot\|_{\infty}, a b^{-l}) \leq&\#(\Omega|_{B_S(m)})^{\lfloor w l\rfloor} \cdot\#(\Omega^{\prime}|_{B_S(m)})^{l-\lfloor wl\rfloor}.
		\end{align*}

		Here $\#(\Omega|_{B_S(m)}) $ and  $\#(\Omega^{\prime}|_{B_S(m)})$ denote the cardinality of  $\Omega|_{ B_S(m)}$  and  $\Omega^{\prime}|_{ B_S(m)} $ respectively.
	\end{lem}
	\begin{proof}
		For any  $\left.v \in \Omega^{\prime}\right|_{B_S(m)}$,  we choose  $s(v) \in A^{B_S(m)} $ satisfying  $\left.(s(v), v) \in \Omega\right|_{B_S(m)} $. Fix $\left.(\xi, \eta) \in \Omega\right|_{B_S(m)} $. for $\left.\left(x_{n}, y_{n}\right) \in \Omega\right|_{B_S(m)}(1\leq n \leq \lfloor wl\rfloor) $ and $\left.y_{n} \in \Omega^{\prime}\right|_{ B_S(m)}(\lfloor wl\rfloor+1\leq n \leq l)$ , write
		$$q((x_{i})_{i=1}^{\lfloor wl\rfloor}, (y_{i})_{i=1}^l)=\left(\sum_{i=1}^{\lfloor wl\rfloor} \frac{x_{n}}{a^{n}}+\sum_{n=\lfloor wl\rfloor+1}^{l} \frac{s\left(y_{n}\right)}{a^{n}}+\sum_{n=l+1}^{\infty} \frac{\xi}{a^{n}}, \sum_{n=1}^l \frac{y_{n}}{b^{n}}+\sum_{n=l+1}^\infty \frac{\eta}{b^{n}}\right).$$
		Then we have $q((x_{i})_{i=1}^{\lfloor wl\rfloor}, (y_{i})_{i=1}^l)\in \left.X_{\Omega}\right|_{ B_S(m)} $. Indeed, if $$((x_{i})_{i=1}^{\lfloor wl\rfloor}, (y_{i})_{i=1}^l)\neq ((x_{i}^\prime)_{i=1}^{\lfloor wl\rfloor}, (y_{i}^\prime)_{i=1}^l),$$ then
		$$b^{-l}=\min \left(a^{-\lfloor w l\rfloor}, b^{-l}\right) \leq \left\|q((x_{i})_{i=1}^{\lfloor wl\rfloor}, (y_{i})_{i=1}^l)-q((x_{i}^\prime)_{i=1}^{\lfloor wl\rfloor}, (y_{i}^\prime)_{i=1}^l)\right\|_{\infty}.$$
		Therefore 
		$$\left\{q((x_{i})_{i=1}^{\lfloor wl\rfloor}, (y_{i})_{i=1}^l) \mid \begin{array}{l}
			\left.\left(x_{m}, y_{m}\right) \in \Omega\right|_{ B_S(m)} \text { for } 1\leq m\leq \lfloor wl\rfloor, \\
			\left.y_{m} \in \Omega^{\prime}\right|_{ B_S(m)} \text { for } \lfloor wl\rfloor\leq m\leq l
		\end{array}\right\}$$
		
		is  $b^{-l}$ -separated set of $\left.X_{\Omega}\right|_{ B_S(m)} $ with respect to the  $\ell^{\infty} $-distance. This gives the first inequality.
		
		\hspace{4mm}
		Since
		$$\bigcup\left\{\begin{array}{l|l}
			\Phi_{l, m}((x_{i})_{i=1}^{\lfloor w l\rfloor}, (y_{i})_{i=1}^l) & \begin{array}{l}
				\left.\left(x_{i}, y_{i}\right) \in \Omega\right|_{ B_S(m)} \text { for }  1\leq i \leq \lfloor wl\rfloor, \\
				\left.y_{i} \in \Omega^{\prime}\right|_{ B_S(m)} \text { for } \lfloor wl\rfloor+1\leq i \leq l
			\end{array}
		\end{array}\right\} $$		
		is an open cover of $X_{\Omega}\mid_{ B_S(m)}$ with 
		$$\operatorname{diam}\left(\Phi_{l, m}(x, y),\|\cdot\|_{\infty}\right) \leq \max \{a^{-\lfloor wl\rfloor}, b^{-l}\}=a^{-\lfloor wl\rfloor}<a b^{-l} .$$
		This gives the second inequality.
		
	\end{proof}

	\begin{pro}
		Let $s$  be a natural number satisfying  $\sum_{g\in G\setminus B_S(s)} \alpha_g<\varepsilon / 2 ,$ and let $c=\sum_{g\in B_S(s)}\alpha_g$. For any $m, l\in \mathbb{N}$, we have
		$$S(X_\Omega,G,d, 2cab^{-l})\leq \frac{\log \left(\#(\Omega|_{B_S(m+s)})^{\lfloor wl\rfloor} \cdot\#(\Omega^{\prime}|_{ B_S(m+s)})^{l-\lfloor w l\rfloor}\right)}{-\log (2ca)+ l\log b}.$$
	\end{pro}
	\begin{proof}
		For any $x, y\in X|_{\Omega}$, and suppose $x=(p_g, q_g)_{g\in G}, y=(p_g^{\prime}, q_g^\prime)_{g\in G},$ then
		
		\begin{align*}
			d_{B_S(m)}(x, y)&=\max_{h\in B_S(m)}d(h((p_g, q_g)_{g\in G}, h(p_g^{\prime}, q_g^\prime)_{g\in G})\\
			&=\max_{h\in B_S(m)} \sum_{g\in G} \alpha_g \max \left(\left|p_{gh}-p_{gh}^{\prime}\right|,\left|q_{gh}-q_{gh}^{\prime}\right|\right)\\
			&=\max_{h\in B_S(m)}\{\sum_{g\in B_S(s)} \alpha_g \max \left(\left|p_{gh}-p_{gh}^{\prime}\right|,\left|q_{gh}-q_{gh}^{\prime}\right|\right)\\
			&\quad +\sum_{g\in G\setminus B_S(s)} \alpha_g \max \left(\left|p_{gh}-p_{gh}^{\prime}\right|,\left|q_{gh}-q_{gh}^{\prime}\right|\right)\}\\
			&\leq \max_{h\in B_S(m+s)}(|p_{h}-p_{h}^{\prime}|,\left|q_{h}-q_{h}^{\prime}\right|)\cdot \sum_{g\in B_{S}(s) } \alpha_g+\sum_{ g\in G\setminus B_S(s)} \alpha_g\\
			&\leq c\left\|\left.x\right|_{B_S(m+s)}-\left.y\right|_{B_S(m+s)}\right\|_{\infty}+\frac{\varepsilon}{2},
		\end{align*}
		thus we have
		
		$$\#\left(X_{\Omega}, d_{B_S(m)}, \varepsilon\right) \leq \#\left(\left.X_{\Omega}\right|_{ B_S(m+s)},\|\cdot\|_{\infty}, \frac{\varepsilon}{2c}\right) .$$

		By Lemma \ref{l5.6}, 
		\begin{align*}
			S(X_\Omega,G,d, 2cab^{-l})&=\frac{\log \#(X_\Omega, d_{B_S(m)}, 2cab^{-l})}{|B_S(m)|}\\
			&\leq \frac{\log \#(X_{\Omega}|_{ B_S(m+s)},\|\cdot\|_{\infty}, ab^{-l})}{|B_S(m)|}\\
			&\leq \frac{\log \left(\#(\Omega|_{B_S(m+s)})^{\lfloor wl\rfloor} \cdot\#(\Omega^{\prime}|_{ B_S(m+s)})^{l-\lfloor w l\rfloor}\right)}{|B_S(m)|}. \qedhere
		\end{align*}
	\end{proof}
	
	\begin{cor}\label{c61}
		Let $R=\sum_{g\in G} \alpha_g$. For any $m,l\in \mathbb{N}$, we have
		$$S(X_\Omega, G,d,2Rab^{-l})\leq \frac{\log \left(\#(\Omega|_{B_S(m+s)})^{\lfloor wl\rfloor} \cdot\#(\Omega^{\prime}|_{ B_S(m+s)})^{l-\lfloor w l\rfloor}\right)}{|B_S(m)|}.$$
	\end{cor}

	\begin{thm}\label{t5.3} 
		The metric mean dimension of the carpet system $\left(X_{\Omega}, G \right)$ is given by
		$$	\operatorname{mdim}_{\mathrm{M}}\left(X_{\Omega}, G, d\right)  =\frac{h_{\rm{top }}(\Omega,G)}{\log a}+\left(\frac{1}{\log b}-\frac{1}{\log a}\right) h_{\rm{top}}\left(\Omega^{\prime}, G\right).$$
	\end{thm}
	
	\begin{proof}
		By Corollary \ref{c61}, 
		\begin{align*}
			&\lim_{m\to \infty}\frac{S(X, G,d,2Rab^{-l})}{|B_S(m)|}\\
			&\leq \lim_{m\to \infty} \frac{\log \left(\#(\Omega|_{ B_S(m+s)})^{\lfloor wl\rfloor} \cdot\#(\Omega^{\prime}|_{ B_S(m+s)})^{l-\lfloor w l\rfloor}\right)}{|B_S(m+s)|}\cdot \frac{| B_S(m+s)|}{| B_S(m)|}\\
			&= \lfloor wl\rfloor h_{\text {top }}(\Omega, G)+(l-\lfloor wl\rfloor) h_{\text {top }}\left(\Omega^{\prime}, G\right)
		\end{align*}
		
		$2Rab^{-l}$ goes to zero as  $L$ goes to infinity, so
		\begin{align*}
			\overline{\operatorname{mdim}}_{\mathrm{M}}\left(X_{\Omega}, G, d\right)&= \limsup_{l \rightarrow \infty} \frac{S(X_\Omega, G, d, 2Rab^{-l})}{-\log (2Ra)+ l\log b}\\
			& \leq\limsup_{l \rightarrow \infty} \frac{\lfloor wl\rfloor h_{\text {top }}(\Omega, G)+(l-\lfloor wl\rfloor) h_{\text {top }}\left(\Omega^{\prime}, G\right)}{-\log( 2Ra)+ l\log b} \\
			& =\frac{w}{\log b} h_{\text {top}}(\Omega, G)+\frac{1-w}{\log b} h_{\text {top}}\left(\Omega^{\prime}, G\right) \\
			& =\frac{h_{\text {top }}(\Omega, G)}{\log a}+\left(\frac{1}{\log b}-\frac{1}{\log a}\right) h_{\text {top }}\left(\Omega^{\prime}, G\right). 
		\end{align*}
		The rest task is to prove the other side.
		Let  $\varepsilon $ be a positive number. For any  $m \geq 1$, since 
		$$\left\|\left.x\right|_{B_S(m)}-\left.y\right|_{B_S(m)}\right\|_{\infty} \leq d_{B_S(m)}(x, y) \quad\left(x, y \in X_{\Omega}\right),$$
		
		where $x|_{B_S(m)}$  denotes the projection of  $x$  to  $\left.X_{\Omega}\right|_{B_S(m)}$. Then for any $\epsilon>0$, we have
		$$\#\left(\left.X_{\Omega}\right|_{B_S(m)},\|\cdot\|_{\infty}, \varepsilon\right) \leq \#\left(X_{\Omega}, d_{B_S(m)}, \varepsilon\right).$$
		By Lemma \ref{l5.6}, we have 
		$$\log \#(X_{\Omega}|_{B_S(m)},\|\cdot\|_{\infty}, b^{-l}) \geq \log (\#(\Omega|_{B_S(m)})^{\lfloor wl\rfloor} \cdot\#(\Omega^{\prime}|_{B_S(m)})^{l-\lfloor w l\rfloor},$$
		thus
		$$\frac{\log \#(X_\Omega, d_{ B_S(m)},  b^{-l})}{|B_S(m)|l\log b}\geq \frac{\lfloor wl\rfloor \log \#(\Omega|_{B_S(m)})+(l-\lfloor wl\rfloor) \log \#( \Omega^{\prime}|_{B_S(m)} )}{|B_S(m)|l\log b}.$$
		Letting $m \rightarrow \infty$,
		$$\frac{S(X, G, d,  b^{-l})}{l\log b}\geq \liminf _{l \rightarrow \infty} \frac{\lfloor wl\rfloor h_{\text {top }}(\Omega, G)+(l-\lfloor wl\rfloor) h_{\text {top }}\left(\Omega^{\prime}, G\right)}{l \log b}.$$
		So
		$$\liminf_{l \rightarrow \infty}\frac{S(X, G, d,  b^{-l})}{l\log b}\geq \frac{w}{\log b} h_{\text {top}}(\Omega, G)+\frac{1-w}{\log b} h_{\text {top}}\left(\Omega^{\prime}, G\right).$$
		This proves \begin{equation*}\underline{\operatorname{mdim}}_{\mathrm{M}}\left(X_{\Omega}, G, d\right)\geq \frac{h_{\text {top }}(\Omega, G)}{\log a}+\left(\frac{1}{\log b}-\frac{1}{\log a}\right) h_{\text {top }}\left(\Omega^{\prime}, G\right).\qedhere
		\end{equation*}
	\end{proof}
	
	Finally, we obtain our conclusion.
	
	\begin{thm}(=Theorem \ref{t1.1})
		For the carpet system $(X_{\Omega},G)$, we have
		\begin{align*}
			&\operatorname{mdim}_{\mathrm{H}}\left(X_{\Omega}, \{B_S(m)\},  d\right) =\frac{h_{\text {top }}^{w}(\pi, G)}{\log b},\\
			&\operatorname{mdim}_{\mathrm{M}}\left(X_{\Omega}, G, d\right)= \frac{h_{\text {top }}(\Omega, G)}{\log a}+\left(\frac{1}{\log b}-\frac{1}{\log a}\right) h_{\text {top }}\left(\Omega^{\prime}, G\right).
		\end{align*}
		Here $\{B_S(m)\}$ is a F$\phi$lner sequence of $G$. 
	\end{thm}
	
	\begin{remark}
		In this section, we failed to study carpet for general amenable group action, but in special finitely generated amenable groups. We know that weighted topological entropy is independent of the choice of metrics, while we calculate the weighted topological entropy of the factor map, an appropriate metric is critical. As in \cite{WCZ}, here we borrow the tool ``length" $l_S(\cdot)$ in finitely generated amenable groups to achieve our purpose.	\end{remark}

	\appendix
	\section{Two examples}
	\hspace{4mm}
	In this appendix, we will give two examples to calculate the metric mean dimension and mean Hausdorff dimension. The following example is where metric mean dimension and mean Hausdorff dimensions are not equal.
	\begin{thm}\label{t1}
		Let $K=\{0\}\cup\{\frac1n\mid n\geq1\}=\{0,1,\frac12,\frac13,\ldots\}$ and $G$ be a countable discrete amenable group. 
		The full $G$-shift $\sigma$ on $
		K^{G}$ is the $G$-action $(K^{G},\sigma)$  defined by
		$$
		\sigma: G \times K^{G} \rightarrow K^{G}
		,(h,(x_{g})_{g \in G}) \mapsto (x_{gh})_{g\in G}, 
		$$
		and the metric $d$ on $K^{G}$ is defined by
		$$d\left((x_g)_{g\in G},(y_g)_{g\in G}\right)=\sum_{g \in G}\alpha_{g}|x_g-y_g|,$$
		where
		$\alpha_{g} \in (0,+\infty)$ satisfies $$
		\alpha_{1_{G}}=1,\sum_{g \in G}\alpha_{g} < +\infty.
		$$
		
		Let $\{F_{n}\}^{\infty}_{n=1}$ be a F$\phi$lner sequence in $G$, then
		\begin{center}
			$\mathrm{mdim}_{\mathrm{H}}\left(K^{G},\{{F_{n}}\},d\right)=0,\quad\mathrm{mdim}_{\mathrm{M}}\left(K^{G},G,d\right)=\frac12.$
		\end{center}
	\end{thm}
	\hspace{4mm}
	To calculate the mean Hausdorff dimension, we give a main technical lemma which is proved by Tsukamoto \cite{TM}.
	\begin{lem}\cite{TM}\label{l13}
		Let $\varepsilon$ and $s$ be positive numbers with $\varepsilon <\frac{1}{6}$. Let $( X, d) $ be a compact metric space with a Borel  probability measure $\mu$,  If for any $x\in X$ there exists a Borel subset $A\subset X$ with $ x\in A$ satisfying
		$$0<\operatorname{diam}A<\frac{\varepsilon}{6},\quad\mu(A) \geq\left(\operatorname{diam}A\right)^s,$$
		then $\dim_{\mathrm{H}}(X,d,\varepsilon)\le2s.$ 
	\end{lem}
	\begin{pro}
		Let $(K^{G},G,\sigma)$ be the  $G$-system in Theorem \ref{t1}, then we have
		$$\overline{\mathrm{mdim}}_{\mathrm{H}} ( K^{G} , \{F_{n}\} , d ) = 0.$$
	\end{pro}
	\begin{proof}
		Set $X=K^{G}$. We define a probability measure $\nu$ on $K$ by
		$$
		\nu\left(\{u\}\right)=au^2\quad(u\neq0),\quad\nu\left(\{0\}\right)=\frac12,
		$$
		where $a$ is a positive number satisfying
		$$a\left(1+\frac1{2^2}+\frac1{3^2}+\cdots\right)=\frac12.$$ 
		
		It is easy to know  $a=\frac{3}{\pi^{2}}<1$. We consider a Borel probability measure on $X$ which is defined by $\mu=\nu^{\otimes G}.$ 
		
		\hspace{4mm}
		Let $c=\sum_{g \in G}\alpha_g>1$. Let  $0 < \varepsilon< \frac 16$ and  $k$ be a natural number. We can find a positive number $\delta=\delta(\varepsilon,k)$ satisfying
		$$\delta<\min\left(\frac\varepsilon{12},\frac{a^k}{(1+c)^{3}},\frac{1}{1+c}\left(\frac12\right)^{\frac k3+1}\right). \quad \quad (1)$$
		Take $S \in \text{Fin}(G)$ satisfying
		$$\sum_{g\in G\setminus S}\alpha_{g}<\delta^{k^k}.$$
		\begin{claim} \label{c14} 
			For any $x\in X$ and any natural munber $n$, there exists a Borel subset $A\subset X$ with $x\in A$ satisfying
			$$
			0<\operatorname{diam}\left(A,d_{F_{n}}\right)<\frac{\varepsilon}{6},\quad\mu(A)\geq\left(\operatorname{diam}(A,d_{F_{n}})\right)^{\frac{6}{k}|SF_{n}|}.
			$$
		\end{claim}
		\begin{proof}
			Take a point $x=(x_g)_{g\in G}$ with $x_g\in K$. For this $x$ we consider a partition
			$$ SF_{n} = I_0\cup I_1\cup I_2\cup \cdots \cup I_k\cup I_{k+ 1}\quad \text{(disjoint union)}$$
			by
			$$I_0=\left\{\left.g\right|x_g>\delta\right\},\quad I_m=\left\{\left.g\right|\delta^{k^m}<x_g\leq\delta^{k^{m-1}}\right\}\quad(1\leq m\leq k),$$
			$$I_{k+1}=\left\{\left.g\right|x_g\leq\delta^{k^k}\right\}.$$
			There exists $k_0\in\{0,1,2,\ldots,k\}$ satisfying
			$$|I_{k_0}|\le\frac{|SF_{n}|}{k+1}<\frac{|SF_{n}|}{k}. \quad \quad (2)$$
			Set
			$$r=\delta^{k^{k_0}}.$$
			We have $r \leq\delta<\varepsilon/(12c).$ We define $A\subset X$ by $$A=\left\{\left.(y_g)_{g\in G}\in K^{G}\right|x_g-r\leq y_g\leq x_g\text{ for all } g \in SF_{n}\right\}.$$
			This is not a single point which implies its diameter is positive. It is obvious that $x\in A$ and
			$$
			\mathrm{diam}(A,d_{F_{n}})\leq cr+\sum_{g\in G\setminus S}\alpha_{g}<cr+\delta^{k^k}\leq (1+c)r < \frac\varepsilon6.
			$$
			Now we  will estimate
			\begin{center}
				$\mu(A) =\prod_{g \in SF_{n}}\nu\left(\left[x_g-r,x_g\right]\right).$ 
			\end{center}
			Let $g \in  SF_{n}$, we need to consider three situations.
			
			Case 1. If $g \in I_{m}$ with $1 \leq m < k_{0}$, then $x_{g} > \d^{k^{m}}$, hence
			\begin{align*}
				\n([x_{g}-r,x_{g}])&\geq \n(\{x_{g}\}) =ax^{2}_{g} \\ &>a\d^{2k^{m}} \\
				&\geq a\d^{2k^{k_{0}-1}}=a(\d^{k^{k_{0}}})^{\frac{2}{k}}=ar^{\frac{2}{k}} \\
				&>((1+c)r)^{\frac{3}{k}} \quad (\text{by }r \leq \d < \frac{a^{k}}{(1+c)^{3}}\text{ in }(1)) \\
				&\geq (\text{diam}(A,d_{F_{n}}))^{\frac{3}{k}}.
			\end{align*}

			Case 2. If $g \in I_{k_{0}}$, then $x_{g} > \d^{k^{k_{0}}}=r$. Hence
			\begin{align*}
				\n([x_{g}-r,x_{g}])&\geq \n(\{x_{g}\})=ax^{2}_{g} \\
				&> ar^{2}\\
				&> ((1+c)r)^{3} \quad (\text{by}~ r \leq \d <\frac{a^{k}}{(1+c)^{3}}< \frac{a}{(1+c)^{3}} )  \\
				&\geq (\text{diam}(A,d_{F_{n}}))^{3}.
			\end{align*}
			Since we know that ${|I_{k_{0}}|< \frac{|SF_{n}|}{k}}$ by $(2)$, we have 
			\begin{center}
				$\prod_{g\in I_{k_{0}}}\n([x_{g}-r,x_{g}])\geq (\text{diam}(A,d_{F_{n}}))^{\frac{3}{k}|SF_{n}|}$.
			\end{center}
			
			Case 3. If $g \in I_{m}$ with $m> k_{0}$, then $x_{g} \leq \d^{k^{m-1}} \leq \d^{k^{k_{0}}}=r$. So $0 \in [x_{g}-r,x_{g}]$. Therefore
			\begin{align*}
				\n([x_{g}-r,x_{g}]) & \geq \n(\{0\})=\frac{1}{2} \\
				&>((1+c)r)^{\frac{3}{k}} \quad ( \text{by } r \leq \d <\frac{1}{1+c} \cdot \left(\frac{1}{2}\right)^{\frac{k}{3}+1}) \\
				& \geq (\text{diam}(A,d_{F_{n}}))^{\frac{3}{k}}.
			\end{align*}
			Summarizing the above, we have
			\begin{align*}
				\m(A) &\geq  (\text{diam}(A,d_{F_{n}}))^{\frac{3}{k}|SF_{n}|} \cdot (\text{diam}(A,d_{F_{n}}))^{\frac{3}{k}|SF_{n}|} \\
				&=(\text{diam}(A,d_{F_{n}}))^{\frac{6}{k}|SF_{n}|}.
			\end{align*}
			This finishes the proof of this claim.
		\end{proof}
		
		\hspace{4mm}
		Combining Lemma \ref{l13} and Claim \ref{c14}, we have
		$$
		\dim_{\rm{H}}(X,d_{F_{n}},\varepsilon)\leq\frac{12}{k}|SF_{n}|
		$$
		for any natural number $n \geq 1$, which implies
		$$\limsup_{n\to\infty}\left(\frac{\dim_{\mathrm{H}}(X,d_{F_{n}},\varepsilon)}{|F_{n}|}\right)\leq 
		\limsup_{n\to\infty} \frac{|SF_{n}|}{|F_{n}|}\frac{12}{k}=\frac{12}{k}.$$
		Letting $\varepsilon$ go to zero and $k$ go to infinity, we get
		\begin{equation*}
			\overline{\mathrm{mdim}}_{\mathrm{H}} ( X ,\{F_{n}\} , d ) = 0.    \qedhere
		\end{equation*}
	\end{proof}

	\begin{thm}
		
		Let $\left(K^{G},G,\sigma\right)$ be the $G$-system in Theorem \ref{t1}, then we have
		$$\quad\mathrm{mdim}_{\mathrm{M}}\left(K^{G},G,d\right)=\frac12.$$
	\end{thm}
	
	\begin{proof}
		Let $0 <\v < \frac{1}{2}$ and $c=\sum_{g \in G}\alpha_{g} < +\infty$. Take $S \in \text{Fin}(G)$ such that 
		$$
		\sum_{g \in G \setminus S}\alpha_{g} \leq \frac{\v}{2}.
		$$
		
		Then we can find a positive integer $\zeta$ which satisfies
		
		$$\frac{1}{\zeta(\zeta+1)} \leq \frac{\v}{4c} < \frac{1}{\zeta(\zeta-1)}.$$

		Clearly, we can use  $2\zeta$ open intervals of $[0,1]$ with length $\frac{\v}{4c}$  to cover $K$. Suppose these open sets are $L_{1},L_{2},...,L_{2\zeta}$. 
		
		\hspace{4mm}
		Let $n$ be any natural number. We consider the following open cover $\mathcal{M}$ of $K^{G}$:
		$$
		\{L_{i_{1}} \cap K\} \times
		\{L_{i_{2}} \cap K\} \times \{L_{i_{3}} \cap K\} \times \cdots\times \{L_{i_{|SF_{n}|}} \cap K\} \times K \times K\cdots,
		$$
		where $i_{j} \in \{1,2,...,2\zeta\}$, $j=\{1,2,...,|SF_{n}|\}$. 
		
		\hspace{4mm}
		Given an open set $M \in \mathcal{M}$. For any  $x=(x_{g})_{g \in G}, y=(y_{g})_{g \in G} \in M$, we have
		\begin{align*}
			&d_{F_{n}}((x_{g})_{g \in G},(y_{g})_{g \in G}) \\
			&=\max_{h \in  F_{n}} d((x_{gh})_{g \in G},(y_{gh})_{g \in G}) \\
			&=\max_{h \in  F_{n} }\sum_{g \in G}\alpha_{g} |x_{gh}-y_{gh}|\\
			&= \max_{h \in  F_{n}} \{\sum_{g \in S}\alpha_{g}{|x_{gh}-y_{gh}|}+\sum_{g \in G \setminus S}\alpha_{g} {|x_{gh}-y_{gh}|}\} \\
			&\leq \sum_{g \in G} \alpha_{g} \frac{\v}{4c}+\frac{\v}{2} \\
			&<\v.
		\end{align*}
		Therefore, for any $M \in \mathcal{M}$, we have $\text{diam}(M,d_{F_{n}}) < \varepsilon$, which implies 
		$$\#(K^{G},d_{F_{n} },\v) \leq (2\zeta)^{|SF_{n}|} .$$
		
		\hspace{4mm}
		Since $\{F_{n}\}^{\infty}_{n=1}$ is a F$\phi$lner sequence, then $\lim_{n \rightarrow \infty}\frac{|SF_{n}|}{|F_{n}|}=1$. Hence
		$$\lim_{n \rightarrow \infty}\frac{\log\#(K^{G},d_{F_{n} },\v)}{|F_{n}|\log {(1/\v)}} \leq \lim_{n \rightarrow \infty}\frac{\log (2\zeta)^{|SF_{n}|}}{|F_{n}|\log{(1/\v)}} 
		= \lim_{n \rightarrow \infty}\frac{
			|SF_{n}|
			\log (2\zeta)}{|F_{n}|\log (1/\v)} 
		=\frac{\log (2\zeta)}{\log (1/ \varepsilon)}.
		$$
		Since $\zeta$ goes to infinity as $\v$ goes to zero, we have
		$$\overline{\text{mdim}}_{\text{M}}(K^{G},G,d) =\limsup_{\v \rightarrow 0}\left(\lim_{n \rightarrow \infty}\frac{\log\#(K^{G},d_{F_{n} },\v)}{|F_{n}|\log {(1/\v)}}\right)\leq\limsup_{\zeta \rightarrow \infty}\frac{\log (2\zeta)}{\log \zeta(1+\zeta)}=\frac{1}{2}.
		$$
		\hspace{4mm}
		On the other hand, given $0 <\v < \frac{1}{4}$, we choose a positive integer $\gamma$ satisfying
		
		$$ \frac{1}{\gamma+1} \leq 2\sqrt{\v}< \frac{1}{\gamma}.$$
		
		\hspace{4mm}
		Fix $z=(z_{g})_{g \in G}$. We consider the sets $W$, where 
		$$W=\{x=(x_g)_{g\in G}\in K^{G}\mid x_g\in \{0,\frac{1}{\gamma},\frac{1}{\gamma-1},\cdots,1\} \text{ for all } g\in F_{n}, \text{ otherwise } x_g=z_{g}\}.$$
		Then for any $x=(x_{g})_{g \in G}, z=(z_{g})_{g \in G}\in W $, we have
		$$ d_{F_{n}}((x_{g})_{g \in G},(z_{g})_{g \in G}) 
		=\max_{h \in  F_{n}} \sum_{g \in G}\alpha_{g}|x_{gh}-z_{gh}| 
		\geq \max_{h\in F_{n} }|x_h-z_h| >\v.
		$$
		Hence, for any distinct points in the sets $W$ have distance $>\varepsilon$ with respect to $d_{F_{n}}$, which implies
		$$\#(K^{G},d_{F_{n}},\v) \geq (\gamma+1)^{|F_{n}|}\geq(\frac{1}{2\sqrt{\v}})^{|F_{n}|} .$$
		
		Therefore
		$$\liminf_{\v \rightarrow 0}\lim_{n \rightarrow \infty}\frac{\log\#(K^{G},d_{F_{n}},\v)}{|F_{n}|\log{(1/\v)}} 
		\geq \lim_{n \rightarrow \infty}\frac{\log (\frac{1}{2\sqrt{\v}})^{|F_{n}|}}{|F_{n}|\log{(1/\v)}} =\frac{\log (\frac{1}{2\sqrt{\v}})}{\log{(1/\v)}} .$$	
		Letting $\varepsilon \rightarrow 0$, we have	
		$$\underline{\mathrm{mdim}}_{\mathrm{M}}(K^{G},G,d) =\liminf_{\v \rightarrow 0}\lim_{n \rightarrow \infty}\frac{\log\#(K^{G},d_{F_{n}},\v)}{|F_{n}|\log{(1/\v)}}  \geq \liminf_{\v \rightarrow 0}\frac{\log (\frac{1}{2\sqrt{\v}})}{\log{(1/\v)}} 
		=\frac{1}{2}.$$
		Hence
		$$\quad\mathrm{mdim}_{\text{M}}\left(K^{G},G,d\right)=\frac12.$$
	\end{proof}

	\hspace{4mm}
	Below we present an example where two dimensions are equal.
	
	\begin{thm}\label{t2}
		Let $G$ be a countable discrete amenable group, and let $[0,1]^{G}$ be an infinite dimensional cube, The shift map $\sigma:G\times [0,1]^{G}\to [0,1]^{G}$ and the metric $d$ on $[0,1]^{G}$ are defined as in Theorem \ref{t1}. For $\{F_{n}\}^{\infty}_{n=1}$ be a F$\phi$lner sequence in $G$,
		$$\mathrm{mdim}\left([0,1]^{G},G\right)=\mathrm{mdim}_{\mathrm{M}}\left([0,1]^{G},G,d\right)=\mathrm{mdim}_{\mathrm{H}}\left([0,1]^{G},\{{F_{n}}\},d\right)=1.$$
	\end{thm}
	\begin{proof}
		We first show that $\overline{\mathrm{mdim}}_{\mathrm{M}}([0,1]^{G},G,d) \leq 1$. Let $\varepsilon >0$ and  $c =\sum_{g \in G} \alpha_{g} < +\infty$ 
		. Take a finite nonempty subset $S$ of $G$ such that 
		$$\sum_{g \in G\setminus S}\alpha_{g} \leq \frac{\varepsilon}{2}.$$
		
		\hspace{4mm}
		Consider the open cover  $I=\{I_0, I_1, \cdots, I_{\lfloor \frac{6c}{\v}\rfloor}\}$ of $[0,1]$, where 
		$$I_{i}=(\frac{(i-1)\v}{6c},\frac{(i+1)\v}{6c}).$$
		\hspace{4mm}
		Let $n$ be any natural number. There exists the following open cover $\mathcal{U}$ of $[0,1]^{G}$ :
		$$
		I_{i_{1}}  \times
		I_{i_{2}}  \times I_{i_{3}}  \times \cdots\times I_{i_{|SF_{n}|}}  \times [0,1] \times [0,1]\cdots,
		$$
		where $i_{j} \in \{0,1,2,...,\lfloor \frac{6c}{\v}\rfloor\}$, $j=\{1,2,..,|SF_{n}|\}$. Then for any open set $U \in \mathcal{U}$, and 	for any $x=(x_{g})_{g \in G}, y=(y_{g})_{g \in G}\in U\in \mathcal{U}$, we have 
		\begin{align*}
			& d_{F_{n}}((x_{g})_{g \in G}, (y_{g})_{g \in G})  \\
			&=\max_{h \in  F_{n} }d((x_{gh})_{g \in G}, (y_{gh})_{g \in G}) \\
			&=\max_{h \in  F_{n} }\sum_{g \in G}\alpha_{g}|x_{gh}-y_{gh}| \\
			&= \max_{h \in  F_{n} }\{\sum_{g \in S}\alpha_{g}|x_{gh}-y_{gh}|+\sum_{g \in G\setminus S}\alpha_{g}|x_{gh}-y_{gh}|\} \\
			& \leq \frac{\v}{3c}\cdot\sum_{g \in G}\alpha_{g} +\sum_{g \in G\setminus S} \alpha_{g}  \\
			& = \frac{\v}{3}+\frac{\v}{2} < \v.
		\end{align*}
		Hence

		$$\#([0,1]^{G},d_{F_{n}},\v) \leq|\mathcal{U}|=(1+\lfloor \frac{6c}{\v}\rfloor)^{|SF_{n}|}.$$
		Since $\{F_{n}\}^{\infty}_{n=1}$ is a F$\phi$lner sequence, then $\lim_{n \rightarrow \infty}\frac{|SF_{n}|}{|F_{n}|}=1$. Hence
		$$\lim_{n \rightarrow \infty}\frac{\log\#([0,1]^{G},d_{F_{n}},\v)}{|F_{n}| \log{(1/\v)}}\leq \lim_{n \rightarrow \infty}\frac{{|SF_{n}|}\log(1+\lfloor \frac{6c}{\v}\rfloor)}{|F_{n}|\log{(1/\v)}} 
		=\frac{\log(1+\lfloor \frac{6c}{\v}\rfloor)}{\log{(1/\v)}} .$$
		
		Letting $\varepsilon \rightarrow 0$, we get 
		$$\overline{\mathrm{mdim}}_{\mathrm{M}}([0,1]^{G},G,d)=\limsup_{\v \rightarrow 0}(\lim_{n \rightarrow \infty}\frac{\log\#([0,1]^{G},d_{F_{n}},\v)}{|F_{n}| \log{(1/\v)}})\leq \limsup_{\v \rightarrow 0}\frac{\log(1+\lfloor \frac{6c}{\v}\rfloor)}{\log{(1/\v)}} =1. $$
		\hspace{4mm}
		By Theorem 10.6.1 in \cite{MCO}, we have $\rm{mdim}{([0,1]^{G},G)}=1$.
		Combining with Theorem \ref{x000}, we get the desired result.
	\end{proof}

	\noindent{\bf Acknowledgments} The authors would like to thank their doctoral supervisor, Siming Tu, who provided much guidance and assistance. The authors are grateful to Rui Yang for pertinent comments and suggestions.

\large{ School of Mathematics(Zhuhai), Sun Yat-sen University,\\
Zhuhai, Guangdong, 519000, P.R. China}	
	
\end{document}